\documentclass[biometrika]{imsart}
\usepackage{amsfonts, amsthm}
\usepackage{amssymb}
\usepackage{amsmath}
\usepackage[numbers]{natbib}
\usepackage{graphicx,psfrag,epsf}
\usepackage{enumerate}
\usepackage{url}
\usepackage{caption}
\usepackage{subcaption}
\usepackage{standalone}
\usepackage{sw20elba}
\usepackage{natbib}

\graphicspath{{figures/}}
\renewcommand{\theequation}{\thesection.\arabic{equation}}
\input{tcilatex}
\begin{document}

\begin{frontmatter}

\title{A Critical Value Function Approach, with an Application to
Persistent Time-Series}
\runtitle{A Critical Value Function Approach}


\begin{aug}
\author{\fnms{Marcelo J.} \snm{Moreira}\ead[label=e1]{moreiramj@fgv.br}},
\and
\author{\fnms{Rafael} \snm{Mour\~{a}o}\ead[label=e2]{rmrodrigues@fgvmail.br}}

\runauthor{Moreira and Mour\~{a}o}

\affiliation{FGV/EPGE}

\address{Escola de P\'{o}s-Gradua\c{c}\~{a}o em Economia e Finan\c{c}as (FGV/EPGE)\\
Getulio Vargas Foundation - 11th floor\\
Praia de Botafogo 190\\
Rio de Janeiro - RJ 22250-040\\
\printead{e1}}

\end{aug}

\begin{abstract}
Researchers often rely on the t-statistic to make inference on parameters in
statistical models. It is common practice to obtain critical values by
simulation techniques. This paper proposes a novel numerical method to
obtain an approximately similar test. This test rejects the null hypothesis
when the test statistic is larger than a critical value function (CVF) of
the data. We illustrate this procedure when regressors are highly
persistent, a case in which commonly-used simulation methods encounter
difficulties controlling size uniformly. Our approach works satisfactorily,
controls size, and yields a test which outperforms the two other known similar tests.
\end{abstract}



\end{frontmatter}

\section{Introduction}

The exact distribution of a test statistic is often unknown to the
researcher. If the statistic is asymptotically pivotal, the quantile based
on the limiting distribution is commonly used as a critical value.
Simulation methods aim to find approximations to such a quantile, including
the Quenouille-Tukey jackknife, the bootstrap method of \citet{Efron79}, and
the subsampling approach of \citet{PolitisRomano94}. These methods are valid
for a variety of models, as proved by \citet{BickelFreedman81}, %
\citet{PolitisRomanoWolf99}, and \citet{RomanoShaikh12}.

For some testing problems, however, the distribution of a statistic can
actually be sensitive to the unknown probability distribution. This
dependence can be related to the curvature of the model, in the sense of
\citet{Efron75,
Efron78}, not vanishing asymptotically. This curvature creates difficulty in
drawing inference on parameters of interest. Usual tests may not control
size uniformly and confidence regions do not necessarily cover the true
parameter at the nominal level. As standard simulation methods may not
circumvent the sensitivity to unknown parameters, there is a need to develop
methods that take this dependence into consideration.

The goal of this paper is to directly tackle the core of the problem,
adjusting for the sensitivity of the test statistic to nuisance parameters.
In practice, our method entails replacing a critical value number by a
critical value function (CVF) of the data. The CVF is a linear combination
of density ratios under the null distribution. The weight choice in this
combination is at the discretion of the researcher. We focus here on
combinations that yield tests which are approximately similar. These weights
can be found by linear programming (LP).

We illustrate the CVF procedure with a simple but important model in which
the regressor may be integrated. \citet{Jeganathan95,Jeganathan97} shows
that the limit of experiments is Locally Asymptotically Brownian Functional
(LABF). The usual test statistics are no longer asymptotically pivotal and
the theory of optimal tests based on Locally Asymptotically Mixed Normal
(LAMN) experiments is no longer applicable. The bootstrap and subsampling
may not provide uniform control in size and confidence levels, as discussed
by \citet{BMMRT91} and \citet{PolitisRomano94}, among others. The similar
t-test based on the CVF has null rejection probabilities close to the
nominal level regardless of the value of the autoregressive parameter. While
this paper does not focus on obtaining optimal procedures, the proposed
similar t-test does have good power. It has correct size and outperforms the
two other tests that are known to be similar:\ the $L_{2}$ test of %
\citet{Wright00} and the uniformly most powerful conditionally unbiased
(UMPCU) test of \citet{JanssonMoreira06}.

This paper is organized as follows. Section \ref{t-stat Sec} introduces the
CVF approach and shows how to find tests which are approximately similar in
finite samples. Section \ref{Model Sec} presents the predictive regressor
model and obtains approximations that are asymptotically valid whether the
disturbances are normal or not. Section \ref{Numerical Sec} provides
numerical simulations to compare the similar t-test with existing procedures
in terms of size and power. Section \ref{Conclusion Sec} concludes. Appendix
A provides a bound based on a discrepancy error useful to implement the CVF.
Appendix B contains the proofs to the theory underlying the CVF. The
supplement presents additional numerical results and remaining proofs.

\section{Inference Based on the t-Statistic\label{t-stat Sec}}

Commonly used tests for $H_{0}:\beta \leq \beta _{0}$ against $H_{1}:\beta
>\beta _{0}$ reject the null hypothesis when some statistic $\psi $ is
larger than a critical value $\kappa _{\alpha }$ (we omit here the
dependence of the statistic $\psi $ on the data $R$ for convenience). The
test is said to have size $\alpha $ when the null rejection probability is
no larger than $\alpha $ for any value of the nuisance parameter $\gamma $:%
\begin{equation}
\sup_{\beta \leq \beta _{0},\gamma \in \Gamma }P_{\beta ,\gamma }\left( \psi
>\kappa _{\alpha }\right) =\alpha .  \label{(Size test)}
\end{equation}%
Finding the critical value $\kappa _{\alpha }$ that satisfies (\ref{(Size
test)}) is a difficult task as the null distribution of the statistic $\psi $
is unknown. In practice, the choice of $\kappa _{\alpha }$ is based on the
limiting distribution of $\psi $ for a large sample size $T$. If $\psi $ is
asymptotically pivotal, the null rejection probability at $\beta _{0}$ is
approximately $\alpha $:
\begin{equation}
\lim_{T\rightarrow \infty }P_{\beta _{0},\gamma }\left( \psi >\kappa
_{\alpha }\right) =\alpha .  \label{(Asymptotic distribution)}
\end{equation}

The t-test rejects the null hypothesis when the t-statistic is larger than a
critical value:%
\begin{equation*}
\psi \left( r\right) >\kappa _{\alpha }.
\end{equation*}%
Simulation methods aim to find approximations to the critical value $\kappa
_{\alpha }$. In most cases, the limit of experiments is locally
asymptotically mixed normal (LAMN) and the t-statistic is asymptotically
normal\footnote{%
More precisely, if the family is locally asymptotically quadratic (LAQ) for
all $\gamma $, then the more restrictive LAMN condition must hold for almost
all $\gamma $ (in the Lebesgue sense); see Proposition 4 of
\citet[p.
77]{LeCamYang00}.}. Hence, numerical methods are usually not too sensitive
to the true unknown nuisance parameter. This critical value number can be
simulated using sample draws. The bootstrap approach works if the asymptotic
critical value does not depend on $\gamma $ within the correct neighborhood
around the true parameter value. In some problems, however, the curvature
does not vanish asymptotically. For example, take the case in which the
regressors are persistent. \citet{CavanaghElliottStock95} show that several
standard methods fail to control size uniformly. This failure is due to the
model curvature not vanishing when the series is nearly integrated. In
particular, the asymptotic quantile of the t-statistic does depend on the
nuisance parameters when the autoregressive coefficient $\gamma $ is near
one.

\subsection{The Critical Value Function\label{CVF SubSec}}

We propose a new method to find a test based on the t-statistic which
controls size uniformly. We focus on finding approximately similar tests
where the null rejection probability is close to $\alpha $. The approach
developed here also holds if we instead want to find tests with correct size
(in a uniform sense) in which the null rejection probability is no larger
than $\alpha $.

Our approach takes into consideration the dependence of the null rejection
probability on $\gamma \in \Gamma $ by means of a critical value function.
The proposed test rejects the null hypothesis when%
\begin{equation}
\psi \left( r\right) >\kappa _{\Gamma _{n},\alpha }(r),
\label{(similar t-test)}
\end{equation}%
where $\kappa _{\Gamma _{n},\alpha }$ is the critical value function so that
the test is similar at $\Gamma _{n}=\left\{ \gamma _{1},...,\gamma
_{n}\right\} $, a set containing $n$ values of the nuisance parameter $%
\gamma $. We choose%
\begin{equation}
\kappa _{\Gamma _{n},\alpha }(r)=\sum_{i=1}^{n}k_{i}\frac{f_{\beta
_{0},\gamma _{i}}(r)}{f_{\nu }\left( r\right) },  \label{(CVF)}
\end{equation}%
where $k=\left( k_{1},k_{2},...,k_{n}\right) $ is a vector of dimension $n$
and $f_{\nu }\left( r\right) $ is a density function. By a standard
separation theorem (e.g., \cite{MoreiraMoreira13}), there exist scalars $%
k_{1},k_{2},...,k_{n}$ so that%
\begin{equation}
P_{\beta _{0},\gamma _{i}}\left( \psi \left( R\right) >\kappa _{\Gamma
_{n},\alpha }(R)\right) =\alpha ,  \label{(boundary)}
\end{equation}%
for $i=1,...,n$. Under some regularity conditions, the scalars $%
k_{1},k_{2},...,k_{n}$ are unique.

The critical value function $\kappa _{\Gamma _{n},\alpha }(r)$ guarantees
that the null rejection probability equals $\alpha $ for $\gamma \in \Gamma
_{n}$. With enough discretization, we can ensure that the null rejection
probability does not differ from $\alpha $ by more than a discrepancy error $%
\varepsilon $.

\bigskip

\begin{theorem}
\label{Uniform Size Thm} Assume that $f_{\beta _{0},\gamma }\left( r\right) $
is continuous in $\gamma $ a.e. in $r$. Consider an array $\gamma _{i,n}$ so
that $\sup_{\gamma \in \Gamma }\min_{\gamma _{i,n}\in \Gamma _{n}}\left\vert
\gamma -\gamma _{i,n}\right\vert \rightarrow 0$ as $n\rightarrow \infty $.
Let $\phi _{\Gamma _{n},\alpha }\left( r\right) $ be the test which rejects
the null if
\begin{equation*}
\psi \left( r\right) >\kappa _{\Gamma _{n},\alpha }(r),
\end{equation*}%
where $\Gamma _{n}=\left\{ \gamma _{1,n},...,\gamma _{n,n}\right\} $. Then,
for any $\epsilon >0$, there exists an $n$ large enough so that
\begin{equation*}
\left\vert E_{\beta _{0},\gamma }\phi _{\Gamma _{n},\alpha }\left( r\right)
-\alpha \right\vert \leq \epsilon ,
\end{equation*}%
for any value of $\gamma \in \Gamma $.
\end{theorem}

\bigskip

By Theorem \ref{Uniform Size Thm}, the size is controlled uniformly if the
discretization is fine enough. In practice, we can find the set $\Gamma _{n}$
through an iterative algorithm. We can start by choosing two endpoints in
the boundary of an interval of interest for $\gamma $. The critical value
function ensures that the null rejection probability for those two points is
exactly $\alpha $, but not for the other points in that interval. We can
proceed by computing the null rejection probabilities for these other
points, either analytically or numerically. If at least one of them differs
from $\alpha $ by more than a discrepancy error, we include the point whose
discrepancy is the largest in the next iteration and start the algorithm all
over again. This approach will be discussed in our numerical example in
Section \ref{Numerical Sec}, where we use Monte Carlo simulations to
estimate the null rejection probabilities.

\subsection{The Baseline Density Function\label{Density SubSec}}

By Le Cam's Third Lemma, the critical value function collapses to a constant
if the statistic is asymptotically pivotal and $f_{\nu }\left( r\right) $ is
a mixture of densities when $\beta =\beta _{0}$. Hence, we focus on%
\begin{equation*}
f_{\nu }\left( r\right) =\int f_{\beta _{0},\gamma }\left( r\right) d\nu
\left( \gamma \right) \text{,}
\end{equation*}%
where $\nu $ is a probability measure which weights different values of $%
\gamma $.

One example is a weight, $\nu ^{\ast }$, which assigns the same mass $1/n$
on a finite number of points $\gamma _{1}<\gamma _{2}<...<\gamma _{n}$:%
\begin{equation*}
\nu ^{\ast }\left( d\gamma \right) =\left\{
\begin{array}{cc}
\frac{1}{n} & \text{if }\gamma =\gamma _{i} \\
0 & \text{otherwise}%
\end{array}%
\right. .
\end{equation*}%
The test then rejects the null hypothesis when%
\begin{equation*}
\psi \left( r\right) >\frac{\sum_{i=1}^{n}k_{i}f_{\beta _{0},\gamma _{i}}(r)%
}{n^{-1}\sum_{j=1}^{n}f_{\beta _{0},\gamma _{j}}(r)}.
\end{equation*}

A second probability measure, $\nu ^{\dagger }$, assigns all mass on a point
$\overline{\gamma }$:%
\begin{equation*}
\nu ^{\dagger }\left( d\gamma \right) =\left\{
\begin{array}{cc}
1 & \text{if }\gamma =\overline{\gamma } \\
0 & \text{otherwise}%
\end{array}%
\right. .
\end{equation*}%
The test rejects the null when%
\begin{equation*}
\psi \left( r\right) >\sum_{i=1}^{n}k_{i}\frac{f_{\beta _{0},\gamma _{i}}(r)%
}{f_{\beta _{0},\overline{\gamma }}(r)}.
\end{equation*}

The choice of $\overline{\gamma }$ for the baseline density $f_{\nu
^{\dagger }}\left( r\right) $ can be motivated by those points for which the
limit of experiments is not LAMN. In our leading example in which the
regressor is persistent, it may be natural to set the autoregressive
parameter $\overline{\gamma }$ equal to one. In practice, however, we
recommend using the density average $f_{\nu ^{\ast }}\left( r\right) $. By
choosing the baseline density as a simple average of the densities $f_{\beta
_{0},\gamma _{i}}(r)$, the CVF is bounded and appears to have better
numerical performance for other null rejection probabilities beyond $\Gamma
=\left\{ \gamma _{1},...,\gamma _{n}\right\} $.

\subsection{Linear Programming\label{LP SubSec}}

We can find the multipliers $k_{l}$, $l=1,...,n$, using a linear programming
approach. The method is simple and requires a sample drawn from only one
probability law. Let $R^{\left( j\right) }$, $j=1,...,J$, be i.i.d. random
variables with positive density $f_{\nu ^{\ast }}\left( r\right) $.

Consider the maximization problem%
\begin{eqnarray}
&&\max_{0\leq \phi (R^{\left( j\right) })\leq 1}\text{ }\frac{1}{J}%
\sum_{j=1}^{J}\phi (R^{\left( j\right) })\frac{\psi (R^{\left( j\right)
}).f_{\nu }(R^{\left( j\right) })}{f_{\nu ^{\ast }}(R^{\left( j\right) })}
\label{(M-estimation)} \\
&&\text{s.t. \ }\frac{1}{J}\sum_{j=1}^{J}\phi (R^{\left( j\right) })\frac{%
f_{\beta _{0},\gamma _{l}}(R^{\left( j\right) })}{f_{\nu ^{\ast }}(R^{\left(
j\right) })}=\alpha \text{, }l=1,...,n\text{.}  \notag
\end{eqnarray}%
The solution to this problem is given by the test in (\ref{(similar t-test)}%
) which approximately satisfies the boundary constraints given in (\ref%
{(boundary)}).

We can rewrite the maximization problem as a standard (primal) linear
programming problem:%
\begin{eqnarray}
&&\max_{m_{J}\in \left[ 0,1\right] ^{J}}\text{ }d^{\prime }m_{J}
\label{(primal LP)} \\
&&\text{s.t. \ }Am_{J}=\alpha .1_{n},  \notag
\end{eqnarray}%
where $d=\left( \frac{\psi (R^{\left( 1\right) }).f_{\nu }(R^{\left(
1\right) })}{f_{\nu ^{\ast }}(R^{\left( 1\right) })},...,\frac{\psi
(R^{\left( J\right) }).f_{\nu }(R^{\left( J\right) })}{f_{\nu ^{\ast
}}(R^{\left( J\right) })}\right) ^{\prime }$ and $m_{J}=(\phi (R^{\left(
1\right) }),...,\phi (R^{\left( J\right) }))^{\prime }$ are vectors in $%
\mathbb{R}^{J}$, the $\left( l,j\right) $-entry of the $n\times J$ matrix $A$
is $f_{\beta _{0},\gamma _{l}}(R^{\left( j\right) })/f_{\nu ^{\ast
}}(R^{\left( j\right) })$, and $1_{n}$ is an $n$-dimensional vector with all
coordinates equal to one. The dual program%
\begin{eqnarray}
&&\min_{k_{J}\in \mathbb{R}^{n}}\text{ }1_{n}^{\prime }k_{J}
\label{(dual LP)} \\
&&\text{s.t. \ }A^{\prime }k_{J}=d  \notag
\end{eqnarray}%
yields an approximation to the vector $k$ associated to the critical value
function given in (\ref{(CVF)}).

The approximation $k_{J}$ to the unknown vector $k$ satisfying (\ref%
{(boundary)}) may, of course, be a poor one. This numerical difficulty can
be circumvented if the LLN holds uniformly as $J\rightarrow \infty $ for the
sample averages present in the maximization problem given by (\ref%
{(M-estimation)}). Showing uniform convergence is essentially equivalent to
proving that the functions, appearing in the objective function and the
boundary constraints of (\ref{(M-estimation)}), belong to the
Glivenko-Cantelli class, as defined in \citet[p.
145]{vanderVaart98}. A sufficient condition for uniform convergence is that
the set $\Gamma $ is compact, the functions $f_{\beta _{0},\gamma
}(r)/f_{\nu ^{\ast }}(r)$ are continuous in $\gamma $ for every $r$, and $%
\psi (r)$ is integrable.

\section{Application: Persistent Regressors\label{Model Sec}}

Consider a simple model with persistent regressors. There is a stochastic
equation%
\begin{equation*}
y_{t}=\mu _{y}+\beta x_{t-1}+\epsilon _{t}^{y},
\end{equation*}%
where the variable $y_{t}$ and the regressor $x_{t}$ are observed, and $%
\epsilon _{t}^{y}$ is a disturbance variable, $t=1,...,T$. This equation is
part of a larger model where the regressor can be correlated with the
unknown disturbance. More specifically, we have%
\begin{equation*}
x_{t}=\gamma x_{t-1}+\epsilon _{t}^{x},
\end{equation*}%
where $x_{0}=0$ and the disturbance $\epsilon _{t}^{x}$ is unobserved and
possibly correlated with $\epsilon _{t}^{y}$. We assume that $\epsilon
_{t}=\left( \epsilon _{t}^{y},\epsilon _{t}^{x}\right) \overset{iid}{\sim }%
N\left( 0,\Sigma \right) $ where
\begin{equation*}
\Sigma =\left[
\begin{array}{cc}
\sigma _{yy} & \sigma _{xy} \\
\sigma _{xy} & \sigma _{xx}%
\end{array}%
\right]
\end{equation*}%
is a positive definite matrix. The theoretical results carry out
asymptotically whether $\Sigma $ is known or not. For simplicity, we assume
for now that $\Sigma $ is known.

Our goal is to assess the predictive power of the past value of $x_{t}$ on
the current value of $y_{t}$. For example, a variable observed at time $t-1$
can be used to forecast another in period $t$. The null hypothesis is $%
H_{0}:\beta \leq 0$ and the alternative is $H_{1}:\beta >0$.

For example, consider the t-statistic%
\begin{equation}
\psi =\frac{\widehat{\beta }-\beta _{0}}{\sigma _{yy}^{1/2}\left[
\sum_{t=1}^{T}\left( x_{t-1}^{\mu }\right) ^{2}\right] ^{-1/2}},
\label{(t-statistic)}
\end{equation}%
where $\widehat{\beta }=\sum_{t=1}^{T}x_{t-1}^{\mu
}y_{t}/\sum_{t=1}^{T}\left( x_{t-1}^{\mu }\right) ^{2}$ (hereinafter, $%
x_{t-1}^{\mu }$ is the demeaned value of $x_{t-1}$) and $\beta _{0}=0$. The
nonstandard error variance uses the fact that $\Sigma $ is known. The
one-sided t-test rejects the null hypothesis when the statistic given by (%
\ref{(t-statistic)}) is larger than the $1-\alpha $ quantile of a standard
normal distribution.

If $\left\vert \gamma \right\vert \,<1$, the t-statistic is asymptotically
normal and $\kappa _{\alpha }$ equals $\Phi ^{-1}\left( 1-\alpha \right) $.
If the convergence in the null rejection probability is uniform and the $%
\sup $ (over $\gamma $) and $\lim $ (in $T$) operators in (\ref{(Size test)}%
) and (\ref{(Asymptotic distribution)}) can be interchanged, the test is
asymptotically similar. That is, the test $\psi >\kappa _{\alpha }$ has null
rejection probability close to $\alpha $ for all values of the nuisance
parameter $\gamma $. The admissible values of the auto-regressive parameter $%
\gamma $ play a fundamental role in the uniform convergence. The
discontinuity in the limiting distribution of the t-statistic $\psi $ at $%
\gamma =1$ shows that the convergence is not uniform even if we know that $%
\left\vert \gamma \right\vert \,<1$. The test based on the t-statistic with $%
\kappa _{\alpha }=\Phi ^{-1}\left( 1-\alpha \right) $ has asymptotic size
substantially different from the size based on the pointwise limiting
distribution of $\psi $. One solution to the size problem is to obtain a
larger critical value. The rejection region $\psi >\kappa _{\alpha }$ would
have asymptotically correct size but would not be similar (and would
consequently be biased). We instead apply the CVF approach to circumvent the
size problems in the presence of highly persistent regressors.

\subsection{Finite-Sample Theory\label{Finite-Sample SubSec}}

We want to test $H_{0}:\beta \leq 0$ against $H_{1}:\beta >0$, where the
nuisance parameters are $\mu _{y}$ and $\gamma $. Standard invariance
arguments can eliminate the parameter $\mu _{y}$. Let $y$ and $x$ be $T$%
-dimensional vectors where the $t$-th entries are $y_{t}$ and $x_{t}$,
respectively. Consider the group of translation transformations on the data%
\begin{equation*}
g\circ \left( y,x\right) =\left( y+g.1_{T},x\right) ,
\end{equation*}%
where $g$ is a scalar and $1_{T}$ is the $T-$dimensional vector whose
entries all equal one. This yields a transformation on the parameter space%
\begin{equation*}
g\circ \left( \beta ,\gamma ,\mu _{y}\right) =\left( \beta ,\gamma ,\mu
_{y}+g\right) .
\end{equation*}%
This group action preserves the parameter of interest $\beta $. Because the
group translation preserves the hypothesis testing problem, it is reasonable
to focus on tests that are invariant to translation transformations on $y$.
The test based on the t-statistic is invariant to these transformations.

Any invariant test can be written as a function of the maximal invariant
statistic. Let $q=\left( q_{1},q_{2}\right) $ be an orthogonal $T\times T$
matrix where the first column is given by $q_{1}=1_{T}/\sqrt{T}$. Algebraic
manipulations show that $q_{2}q_{2}^{\prime }=M_{1_{T}}$, where $%
M_{1_{T}}=I_{T}-1_{T}\left( 1_{T}^{\prime }1_{T}\right) ^{-1}1_{T}^{\prime }$
is the projection matrix to the space orthogonal to $1_{T}$. Let $x_{-1}$ be
the $T$-dimensional vector whose $t$-th entry is $x_{t-1}$, and define $%
w^{\mu }=q_{2}^{\prime }w$ for a $T$-dimensional vector $w$. The maximal
invariant statistic is given by $r=\left( y^{\mu },x\right) $. Its density
function is given by
\begin{eqnarray}
&&f_{\beta ,\gamma }\left( y^{\mu },x\right) =\left( 2\pi \sigma
_{xx}\right) ^{-\frac{T}{2}}\exp \left\{ -\frac{1}{2\sigma _{xx}}%
\sum_{t=1}^{T}\left( x_{t}-x_{t-1}\gamma \right) ^{2}\right\}
\label{(likelihood eq)} \\
&&\times \left( 2\pi \sigma _{yy.x}\right) ^{-\frac{T-1}{2}}\exp \left\{ -%
\frac{1}{2\sigma _{yy.x}}\sum_{t=1}^{T}\left( y_{t}^{\mu }-x_{t}^{\mu }\frac{%
\sigma _{xy}}{\sigma _{xx}}-x_{t-1}^{\mu }\left[ \beta -\gamma \frac{\sigma
_{xy}}{\sigma _{xx}}\right] \right) ^{2}\right\} ,  \notag
\end{eqnarray}%
where $\sigma _{yy.x}=\sigma _{yy}-\sigma _{xy}^{2}/\sigma _{xx}$ is the
variance of $\epsilon _{t}^{y}$ not explained by $\epsilon _{t}^{x}$.

Lemma 1 of \citet{JanssonMoreira06} shows that this density belongs to the
curved-exponential family, which consists of two parameters, $\beta $ and $%
\gamma $, and a four-dimensional sufficient statistic:%
\begin{eqnarray*}
S_{\beta } &\hspace{-0.08in}=\hspace{-0.08in}&\frac{1}{\sigma _{yy.x}}%
\sum_{t=1}^{T}x_{t-1}^{\mu }\left( y_{t}-\frac{\sigma _{xy}}{\sigma _{xx}}%
x_{t}\right) \text{, }S_{\gamma }=\frac{1}{\sigma _{xx}}%
\sum_{t=1}^{T}x_{t-1}x_{t}-\frac{\sigma _{xy}}{\sigma _{xx}}S_{\beta }\text{,%
} \\
S_{\beta \beta } &\hspace{-0.08in}=\hspace{-0.08in}&\frac{1}{\sigma _{yy.x}}%
\sum_{t=1}^{T}\left( x_{t-1}^{\mu }\right) ^{2}\text{, and }S_{\gamma \gamma
}=\frac{1}{\sigma _{xx}}\sum_{t=1}^{T}\left( x_{t-1}\right) ^{2}+\frac{%
\sigma _{xy}^{2}}{\sigma _{xx}^{2}}S_{\beta \beta }.
\end{eqnarray*}%
As a result, classical exponential-family statistical theory is not
applicable to this problem. This has several consequences for testing $%
H_{0}:\beta \leq 0$ against $H_{1}:\beta >0$. First, there does not exist a
uniformly most powerful unbiased (UMPU) test for $H_{0}:\beta \leq 0$
against $H_{1}:\beta >0$. \citet{JanssonMoreira06} instead obtain an optimal
test within the class of similar tests conditional on the specific ancillary
statistics $S_{\beta \beta }$ and $S_{\gamma \gamma }$. Second, the pair of
sufficient statistics under the null hypothesis, $S_{\gamma }$ and $%
S_{\gamma \gamma }$, is not complete. Hence, there exist similar tests
beyond those considered by \citet{JanssonMoreira06}. This less restrictive
requirement can lead to power gains because (i) $S_{\beta \beta }$ and $%
S_{\gamma \gamma }$ are not independent of $S_{\beta }$ and $S_{\gamma }$;
and (ii) the joint distribution of $S_{\beta }$ and $S_{\gamma }$ does
depend on the parameter $\beta $.

\subsection{Asymptotic Theory\label{Asymptotic SubSec}}

We now analyze the behavior of testing procedures as $T\rightarrow \infty $.
The asymptotic power of tests depends on the likelihood ratio around $\beta
_{0}=0$ and the true parameter $\gamma $.

The Hessian matrix of the log-likelihood captures the difficulty in making
inference in parametric models. Define the scaling function%
\begin{equation*}
g_{T}\left( \gamma \right) =\left( \frac{1}{\sigma _{xx}}E_{0,\gamma
}\sum_{t=1}^{T}\left( x_{t-1}\right) ^{2}\right) ^{-1/2},
\end{equation*}%
which is directly related to the Hessian matrix of the log-likelihood given
in (\ref{(likelihood eq)}). For example, \citet{Anderson59} and %
\citet{Mikusheva07} implicitly use the scaling function $g_{T}\left( \gamma
\right) $ for the one-equation autoregressive model with one lag, AR(1).
This scaling function is chosen so that $P_{0,\gamma }$ and $P_{b\cdot
\sigma _{yy.x}^{1/2}\sigma _{xx}^{-1/2}g_{T}\left( \gamma \right) ,\gamma
+c\cdot g_{T}\left( \gamma \right) }$ are mutually contiguous for any value
of $\gamma $. The local alternative is indexed by $\sigma
_{yy.x}^{1/2}\sigma _{xx}^{-1/2}$ to simplify algebraic calculations. Define
the stochastic process $\Lambda _{T}$ as the log-likelihood ratio%
\begin{equation*}
\Lambda _{T}\left( b\cdot \sigma _{yy.x}^{1/2}\sigma _{xx}^{-1/2}g_{T}\left(
\gamma \right) ,\gamma +c\cdot g_{T}\left( \gamma \right) ;0,\gamma \right)
=\ln \frac{f_{b\cdot \sigma _{yy.x}^{1/2}\sigma _{xx}^{-1/2}g_{T}\left(
\gamma \right) ,\gamma +c\cdot g_{T}\left( \gamma \right) }\left( r\right) }{%
f_{0,\gamma }\left( r\right) }.
\end{equation*}%
The subscript $T$ indicates the dependency of the log-likelihood ratio on
the sample.

\bigskip

\begin{proposition}
\label{Scaling Fn Prop} The scaling function equals%
\begin{equation*}
g_{T}\left( \gamma \right) =\left( \sum_{t=1}^{T-1}\sum_{l=0}^{t-1}\gamma
{}^{2l}\right) ^{-1/2}.
\end{equation*}%
The function $g_{T}\left( \gamma \right) $ is continuous in $\gamma $ and
simplifies to the usual rates for the stationary, integrated, and explosive
cases:\newline
\emph{(a)} if $\left\vert \gamma \right\vert <1$, then $g_{T}\left( \gamma
\right) \sim (1-\gamma ^{^{2}})^{1/2}T^{-1/2}$;\newline
\emph{(b)} if $\gamma =1$, then $g_{T}\left( \gamma \right) \sim
2^{1/2}T^{-1}$; and,\newline
\emph{(c)} if $\gamma >1$, then $g_{T}\left( \gamma \right) \sim (1-\gamma
^{^{-2}})\cdot \gamma ^{^{-(T-2)}}$.
\end{proposition}

\bigskip

The log-likelihood ratio can be written as%
\begin{equation*}
\Lambda _{T}\left( b\cdot \sigma _{yy.x}^{1/2}\sigma _{xx}^{-1/2}g_{T}\left(
\gamma \right) ,\gamma +c\cdot g_{T}\left( \gamma \right) ;0,\gamma \right)
= \left[ b,c\right] R_{T}\left( \gamma \right) -\frac{1}{2}\left[ b,c\right]
K_{T}\left( \gamma \right) \left[ b,c\right] ^{\prime },
\end{equation*}%
where $R_{T}$ is a random vector and $K_{T}$ is a symmetric random matrix.
The first and second components of $R_{T}$ are%
\begin{eqnarray*}
R_{\beta ,T} &\hspace{-0.08in}\left( \gamma \right) =\hspace{-0.08in}%
&g_{T}\left( \gamma \right) \frac{1}{\sigma _{yy.x}^{1/2}\sigma _{xx}^{1/2}}%
\sum_{t=1}^{T}x_{t-1}^{\mu }\left( y_{t}-\beta _{0}x_{t-1}-\frac{\sigma _{xy}%
}{\sigma _{xx}}\left[ x_{t}-\gamma x_{t-1}\right] \right) \text{ and} \\
R_{\gamma ,T} &\hspace{-0.08in}\left( \gamma \right) =\hspace{-0.08in}%
&g_{T}\left( \gamma \right) \frac{1}{\sigma _{xx}}\sum_{t=1}^{T}x_{t-1}%
\left( x_{t}-\gamma x_{t-1}\right) -\frac{\rho }{\sqrt{1-\rho ^{2}}}R_{\beta
,T}\left( \gamma \right) \text{,}
\end{eqnarray*}%
where $\rho =\sigma _{xy}/(\sigma _{xx}^{1/2}\sigma _{yy}^{1/2})$ is the
error correlation. The entries $(1,1)$, $(1,2)$, and $(2,2)$ of $K_{T}\left(
\gamma \right) $ are, respectively,%
\begin{eqnarray*}
K_{\beta \beta ,T} &\hspace{-0.08in}\left( \gamma \right) =\hspace{-0.08in}%
&g_{T}\left( \gamma \right) ^{2}\frac{1}{\sigma _{xx}}%
\sum_{t=1}^{T}x_{t-1}^{\mu ^{2}}\text{,} \\
K_{\beta \gamma ,T} &\hspace{-0.08in}\left( \gamma \right) =\hspace{-0.08in}%
&-g_{T}\left( \gamma \right) ^{2}\frac{\rho }{\sqrt{1-\rho ^{2}}}\frac{1}{%
\sigma _{xx}}\sum_{t=1}^{T}x_{t-1}^{\mu ^{2}}\text{, and} \\
K_{\gamma \gamma ,T} &\hspace{-0.08in}\left( \gamma \right) =\hspace{-0.08in}%
&g_{T}\left( \gamma \right) ^{2}\left[ \frac{\rho ^{2}}{1-\rho ^{2}}\frac{1}{%
\sigma _{xx}}\sum_{t=1}^{T}x_{t-1}^{\mu ^{2}}+\frac{1}{\sigma _{xx}}%
\sum_{t=1}^{T}x_{t-1}^{^{2}}\right] \text{.}
\end{eqnarray*}%
The asymptotic behavior of the log-likelihood ratio is divided into
pointwise limits: $\left\vert \gamma \right\vert <1$ (stationary), $\gamma
=1 $ (integrated), and $\gamma >1$ (explosive).

\bigskip

\begin{proposition}
\label{Pointwise Approx Prop} For any bounded scalars $b$ and $c$, $\Lambda
_{T}^{\gamma }\left( b,c\right) =\Lambda ^{\gamma }\left( b,c\right)
+o_{P_{0,\gamma }}\left( 1\right) $, where $\Lambda ^{\gamma }\left(
b,c\right) $ is defined as follows:\newline
\emph{(a)} For $\left\vert \gamma \right\vert <1$,%
\begin{equation*}
\Lambda ^{\gamma }\left( b,c\right) =\left[ b,c\right] R^{S}-\frac{1}{2}%
\left[ b,c\right] K^{S}\left[ b,c\right] ^{\prime },
\end{equation*}%
where $R^{S}\sim N\left( 0,K^{S}\right) $ and $K^{S}$ is a matrix given by%
\begin{equation*}
K^{S}=\left[
\begin{array}{cc}
1 & -\frac{\rho }{\sqrt{1-\rho ^{2}}} \\
-\frac{\rho }{\sqrt{1-\rho ^{2}}} & \frac{1}{1-\rho ^{2}}%
\end{array}%
\right] .
\end{equation*}%
\newline
\emph{(b) }For $\gamma =1$,
\begin{equation*}
\Lambda ^{\gamma }\left( b,c\right) =\left[ b,c\right] R^{I}-\frac{1}{2}%
\left[ b,c\right] K^{I}\left[ b,c\right] ^{\prime },
\end{equation*}%
where $R^{I}$ and $K^{I}$ are given by%
\begin{eqnarray*}
R^{I} &\hspace{-0.08in}=\hspace{-0.08in}&2^{1/2}\left[
\begin{array}{c}
\int_{0}^{1}W_{x}^{\mu }\left( r\right) dW_{y}\left( r\right) \\
\int_{0}^{1}W_{x}^{\mu }\left( r\right) dW_{x}\left( r\right) -\frac{\rho }{%
\sqrt{1-\rho ^{2}}}\int_{0}^{1}W_{x}^{\mu }\left( r\right) dW_{y}\left(
r\right)%
\end{array}%
\right] \text{ and} \\
K^{I} &\hspace{-0.08in}=\hspace{-0.08in}&2\left[
\begin{array}{cc}
\int_{0}^{1}W_{x}^{\mu }\left( r\right) ^{2}dr & -\frac{\rho }{\sqrt{1-\rho
^{2}}}\int_{0}^{1}W_{x}^{\mu }\left( r\right) ^{2}dr \\
-\frac{\rho }{\sqrt{1-\rho ^{2}}}\int_{0}^{1}W_{x}^{\mu }\left( r\right)
^{2}dr & \frac{\rho }{1-\rho ^{2}}\int_{0}^{1}W_{x}^{\mu }\left( r\right)
^{2}dr+\int_{0}^{1}W_{x}^{\mu }\left( r\right) ^{2}dr%
\end{array}%
\right] ,
\end{eqnarray*}%
\newline
with $W_{x}$ and $W_{y}$ being two independent Wiener processes and%
\begin{equation*}
W_{x}^{\mu }\left( r\right) =W_{x}\left( r\right) -\int_{0}^{1}W_{x}\left(
s\right) ds.
\end{equation*}%
\newline
\emph{(c) }For $\gamma >1$,
\begin{equation*}
\Lambda ^{\gamma }\left( b,c\right) =\left[ b,c\right] R^{E}-\frac{1}{2}%
\left[ b,c\right] K^{E}\left[ b,c\right] ^{\prime },
\end{equation*}%
where $\left( K^{E}\right) ^{-1/2}R^{E}\sim N\left( 0,I_{2}\right) $
independent of $K^{E}$ whose distribution is%
\begin{equation*}
K^{E}\sim \chi ^{2}\left( 1\right) \cdot \left[
\begin{array}{cc}
1 & -\frac{\rho }{\sqrt{1-\rho ^{2}}} \\
-\frac{\rho }{\sqrt{1-\rho ^{2}}} & \frac{1}{1-\rho ^{2}}%
\end{array}%
\right] .
\end{equation*}
\end{proposition}

The log-likelihood ratio has an exact quadratic form of $b$ and $c$ in
finite samples. Proposition \ref{Pointwise Approx Prop} finds the limit of
experiments when the series is (a) stationary; (b) nearly integrated; and
(c) explosive. When $\left\vert \gamma \right\vert <1$, the limit of
experiments is Locally Asymptotically Normal (LAN). When $\gamma =1$ the
limit of experiments is Locally Asymptotically Brownian Functional (LABF).
When $\gamma >1$, the limit of experiments is Locally Asymptotically Mixed
Normal (LAMN). Because LAN is a special case of LAMN, the limit of
experiments is LAMN\ for every $\gamma $ except for $\gamma =1$. By Le Cam's
Third Lemma, the asymptotic power of a test $\overline{\phi }$ is given by%
\begin{equation*}
\lim_{T\rightarrow \infty }E_{b\cdot \sigma _{yy.x}^{1/2}\sigma
_{xx}^{-1/2}g_{T}\left( \gamma \right) ,\gamma +c\cdot g_{T}\left( \gamma
\right) }\overline{\phi }\left( R\right) =E_{0,\gamma }\overline{\phi }\exp
\left\{ \Lambda ^{\gamma }\left( b,c\right) \right\} .
\end{equation*}

As $T\rightarrow \infty $, the null rejection probability is controlled in a
neighborhood that shrinks with rate $g_{T}\left( \gamma \right) $. Define
\begin{equation*}
f_{0,\gamma _{i}}\left( r\right) =f_{0,\gamma +c_{i}\cdot g_{T}\left( \gamma
\right) }\left( r\right) ,
\end{equation*}%
where $\left\vert c_{i}\right\vert \leq c$ for $i=1,...,n$ (we assume that $%
n $ is odd and $c_{\left( n+1\right) /2}=0$ for convenience). Theorem \ref%
{Uniform Size Thm} and Proposition \ref{Pointwise Approx Prop} allow us to
refine the partition within a shrinking interval for $\gamma $. The
approximation error converges in probability to zero for any bounded values
of $b$ and $c$.

Let $\phi _{\Gamma _{n},\alpha }^{\ast }$ be the test based on the measure $%
\nu ^{\ast }$, which assigns the same mass $1/n$ on $\Gamma _{n}=\left\{
\gamma _{1},...,\gamma _{n}\right\} $. This test rejects $H_{0}$ when%
\begin{equation*}
\psi \left( r\right) >\frac{\sum_{i=1}^{n}k_{i}^{\ast }f_{0,\gamma
_{i}}\left( r\right) }{n^{-1}\sum_{j=1}^{n}f_{0,\gamma _{j}}\left( r\right) }%
,
\end{equation*}%
where $k_{1}^{\ast },...,k_{n}^{\ast }$ are constants such that the null
rejection probability equals $\alpha $ at $\gamma _{1},...,\gamma _{n}$.
Analogously, let $\phi _{\Gamma _{n},\alpha }^{\dagger }$ be the test based
on the measure $\nu ^{\dagger }$, which assigns all mass on $\gamma
_{(n+1)/2}=\gamma $. This test rejects $H_{0}$ when%
\begin{equation*}
\psi \left( r\right) >\frac{\sum_{i=1}^{n}k_{i}^{\dagger }f_{0,\gamma
_{i}}\left( r\right) }{f_{0,\gamma }\left( r\right) },
\end{equation*}%
where $k_{1}^{\ast },...,k_{n}^{\ast }$ are constants such that the null
rejection probability equals $\alpha $ at $\gamma _{i}=\gamma +c_{i}\cdot
g_{T}\left( \gamma \right) $, $i=1,...,n$. We, of course, do not know the
true parameter $\gamma $.

Both tests $\phi _{\Gamma _{n},\alpha }^{\ast }$ and $\phi _{\Gamma
_{n},\alpha }^{\dagger }$ reject $H_{0}$ when the test statistic $\psi
\left( r\right) $ is larger than a critical value function. The respective
critical value functions can be written in terms of the log-likelihood ratio
differences $\Lambda _{T}\left( 0,\gamma +c_{i}\cdot g_{T}\left( \gamma
\right) ;0,\gamma \right) $. The test $\phi _{\Gamma _{n},\alpha }^{\ast }$
rejects $H_{0}$ when%
\begin{equation}
\psi \left( r\right) >\frac{\sum_{i=1}^{n}k_{i}^{\ast }\exp \left\{ \Lambda
_{T}\left( 0,\gamma +c_{i}\cdot g_{T}\left( \gamma \right) ;0,\gamma \right)
\right\} }{n^{-1}\sum_{j=1}^{n}\exp \left\{ \Lambda _{T}\left( 0,\gamma
+c_{j}\cdot g_{T}\left( \gamma \right) ;0,\gamma \right) \right\} }.
\label{(multipliers measure 1)}
\end{equation}%
The test $\phi _{\Gamma _{n},\alpha }^{\dagger }$ rejects $H_{0}$ when%
\begin{equation}
\psi \left( r\right) >\sum_{i=1}^{n}k_{i}^{\dagger }\exp \left\{ \Lambda
_{T}\left( 0,\gamma +c_{i}\cdot g_{T}\left( \gamma \right) ;0,\gamma \right)
\right\} .  \label{(multipliers measure 2)}
\end{equation}

Assume that $\psi \left( R\right) $ is regular and asymptotically standard
normal when the series is stationary (or explosive). In particular, this
assumption applies to the t-statistic. For each $\gamma $, we write $%
\rightarrow _{c}$ to denote convergence in distribution under $P_{0,\gamma
+c.g_{T}\left( \gamma \right) }$.

\bigskip

\noindent \textbf{Assumption AP:} The test statistic $\psi \left( r\right) $
is a continuous function of the sufficient statistics (a.e.). Furthermore,
when $\left\vert \gamma \right\vert <1$ (or $\gamma >1$), $\psi \left(
R\right) \rightarrow _{c}N\left( 0,1\right) $ for any bounded $c$.

\bigskip

Under Assumption AP, the critical value function converges to the $1-\alpha $
quantile of a normal distribution when either $\left\vert \gamma \right\vert
<1$ or $\gamma >1$. This implies that both tests $\phi _{\Gamma _{n},\alpha
}^{\ast }$ and $\phi _{\Gamma _{n},\alpha }^{\dagger }$ are asymptotically
UMPU in the stationary case and conditionally (on ancillary statistics) UMPU
in the explosive case.

\bigskip

\begin{theorem}
\label{Asymptotic Thm} Assume that $\psi $ satisfies Assumption AP. The
following hold as $T\rightarrow \infty $:\newline
\emph{(a)} If $\left\vert \gamma \right\vert <1$, then \emph{(i)} $%
k_{i}^{\ast }\rightarrow _{p}\Phi ^{-1}\left( 1-\alpha \right) /n$, and
\emph{(ii)} $k_{(n+1)/2}^{\dagger }\rightarrow _{p}\Phi ^{-1}\left( 1-\alpha
\right) $ and $k_{i}^{\dagger }\rightarrow _{p}0$ for $i\neq (n+1)/2$.%
\newline
\emph{(b)} If $\gamma =1$, then \emph{(i) }$k_{i}^{\ast }$ converge to $%
k_{i,\infty }^{\ast }$ which satisfy%
\begin{equation}
P\left( \psi \left( c_{l}\right) >\frac{\sum_{i=1}^{n}k_{i,\infty }^{\ast
}\exp \left\{ c_{i}R_{\gamma }\left( c_{l}\right) -\frac{1}{2}%
c_{i}^{2}K_{\gamma \gamma }\left( c_{l}\right) \right\} }{%
n^{-1}\sum_{j=1}^{n}\exp \left\{ c_{j}R_{\gamma }\left( c_{l}\right) -\frac{1%
}{2}c_{j}^{2}K_{\gamma \gamma }\left( c_{l}\right) \right\} }\right) =\alpha
,  \label{(Asymptotically Similar Integrated eq)}
\end{equation}%
$l=1,...,n$, and \emph{(ii)} $k_{i}^{\dagger }$ converge to $k_{i,\infty
}^{\dagger }$ which satisfy%
\begin{equation}
P\left( \psi \left( c_{l}\right) >\sum_{i=1}^{n}k_{i,\infty }^{\dagger }\exp
\left\{ c_{i}R_{\gamma }\left( c_{l}\right) -\frac{1}{2}c_{i}^{2}K_{\gamma
\gamma }\left( c_{l}\right) \right\} \right) =\alpha ,
\label{(Asymptotically Similar Integrated 2 eq)}
\end{equation}%
$l=1,...,n$, where $\psi \left( c_{l}\right) $ is defined as $\psi \left(
R_{\beta }\left( c_{l}\right) ,R_{\gamma }\left( c_{l}\right) ,K_{\beta
\beta }\left( c_{l}\right) ,K_{\beta \gamma }\left( c_{l}\right) ,K_{\gamma
\gamma }\left( c_{l}\right) \right) $ and%
\begin{eqnarray*}
R_{\beta }\left( c\right)  &\hspace{-0.08in}=\hspace{-0.08in}&2^{1/2}\left[
\int_{0}^{1}W_{x,c}^{\mu }\left( r\right) dW_{y}\left( r\right) -\frac{\rho
}{\sqrt{1-\rho ^{2}}}c\int_{0}^{1}W_{x,c}^{\mu }\left( r\right) ^{2}dr\right]
, \\
R_{\gamma }\left( c\right)  &\hspace{-0.08in}=\hspace{-0.08in}&2^{1/2}\left[
\int_{0}^{1}W_{x}^{\mu }\left( r\right) dW_{x}\left( r\right) -\frac{\rho }{%
\sqrt{1-\rho ^{2}}}R_{\beta }\left( c\right) \right] , \\
K_{\beta \beta }\left( c\right)  &\hspace{-0.08in}=\hspace{-0.08in}%
&2\int_{0}^{1}W_{x,c}^{\mu }\left( r\right) ^{2}dr, \\
K_{\beta \gamma }\left( c\right)  &\hspace{-0.08in}=\hspace{-0.08in}&-2\frac{%
\rho }{\sqrt{1-\rho ^{2}}}\int_{0}^{1}W_{x,c}^{\mu }\left( r\right) ^{2}dr,%
\text{ and} \\
K_{\gamma \gamma }\left( c\right)  &\hspace{-0.08in}=\hspace{-0.08in}&2\left[
\frac{\rho }{1-\rho ^{2}}\int_{0}^{1}W_{x,c}^{\mu }\left( r\right)
^{2}dr+\int_{0}^{1}W_{x,c}^{\mu }\left( r\right) ^{2}dr\right] ,
\end{eqnarray*}%
\newline
for the two independent Wiener processes $W_{x}$ and $W_{y}$, with $W_{x,c}$ being the
standard Ornstein-Uhlenbeck process, and $W_{x,c}^{\mu
}=W_{x,c}-\int_{0}^{1}W_{x,c}\left( s\right) ds$.\newline
\emph{(c)} If $\gamma >1$ and $\psi $ is the t-statistic, then \emph{(i)} $%
k_{i}^{\ast }\rightarrow _{p}\Phi ^{-1}\left( 1-\alpha \right) /n$, and
\emph{(ii)} $k_{(n+1)/2}^{\dagger }\rightarrow _{p}\Phi ^{-1}\left( 1-\alpha
\right) $ and $k_{i}^{\dagger }\rightarrow _{p}0$ for $i\neq (n+1)/2$.
\end{theorem}

\textbf{Comment:} Theorem \ref{Asymptotic Thm} is valid for $\left\vert
\gamma \right\vert \leq 1$ when errors are nonnormal and the error
covariance is estimated under mild moment conditions.

\bigskip

Parts (a) and (c) show that there is no critical-value adjustment when the
series is either stationary or explosive. Part (b) provides a size
correction when the series is integrated. We note that the t-statistic is
not asymptotically pivotal under the nearly-integrated asymptotics. Hence,
the size adjustment is non-trivial; that is, not all $n$ boundary conditions
are satisfied for the t-statistic if $k_{i}^{\ast }$ equals a constant for $%
i=1,...,n$ or $k_{i}^{\dagger }=0$ for $i\neq \left( n+1\right) /2$.

We, of course, do not know the true parameter $\gamma $. We could choose an
auxiliary estimator $\widehat{\gamma }_{T}$ or select a prespecified value
of $\gamma $ for the CVF. In this paper, we choose the second option. As the
CVF provides an adjustment only when $\gamma =1$, we recommend this value as
the centering point of the CVF. We use the test
\begin{equation}
\psi \left( r\right) >\frac{\sum_{i=1}^{n}k_{i}^{\ast }\exp \left\{ \Lambda
_{T}\left( 0,1+c_{i}\cdot g_{T}\left( 1\right) ;0,1\right) \right\} }{%
n^{-1}\sum_{j=1}^{n}\exp \left\{ \Lambda _{T}\left( 0,1+c_{j}\cdot
g_{T}\left( \gamma \right) ;0,1\right) \right\} }.  \label{(feasible test)}
\end{equation}%
for the numerical simulations.

\section{Numerical Simulations\label{Numerical Sec}}

In this section, we present numerical simulations for the approximately
similar t-test. We simulate the simple model with moderately small samples ($%
T=100$). We perform 10,000 Monte Carlo simulations to evaluate each
rejection probability. The nominal size $\alpha $ is 10\% and $\beta _{0}=0$%
. The disturbances $\varepsilon _{t}^{y}$ and $\varepsilon _{t}^{x}$ are
serially independent and identically distributed, with variance one and
correlation $\rho $. We focus on the iterative algorithm which provides size
adjustments and on properties of the CVF. To avoid an additional level of
uncertainty, we assume for now the covariance matrix to be known. In Section %
\ref{Power SubSec}, we relax the assumption of known variance and evaluate
the feasible similar t-test in terms of size and power.

To find the approximated critical value function, we start with two initial
points at $c=-50$ and $c=20$. For $T=100$, the value of the autoregressive
parameter is respectively $\gamma =0.5$ and $\gamma =1.2$. The choice of
these parameters relies on the fact that the t-statistic is approximately
normal in the stationary and explosive cases. We uniformly discretize the
range $c\in \left[ -50,20\right] $ by 100 points and make sure all null
rejection probabilities are near 10\%. We approximate the null rejection
probabilities for these points with $J=10,000$ Monte Carlo simulations. If
at least one of them differs from 10\% by more than 0.015, we include the
point whose discrepancy is the largest in the next iteration. We then
re-evaluate all null rejection probabilities and include an additional point
if necessary, and so forth. In Appendix A, we show that if the discrepancy
error is 0.015, the probability of unnecessarily adding one more point is
less than 0.01\%. The multipliers are obtained solving the linear
programming problem given in expression (\ref{(dual LP)}).

\begin{figure}[tbh]
\caption{Null Rejection Probabilities}
\label{fig:Size}
\bigskip \minipage{0.5\textwidth} \centering %
\includegraphics[width=6.8cm]{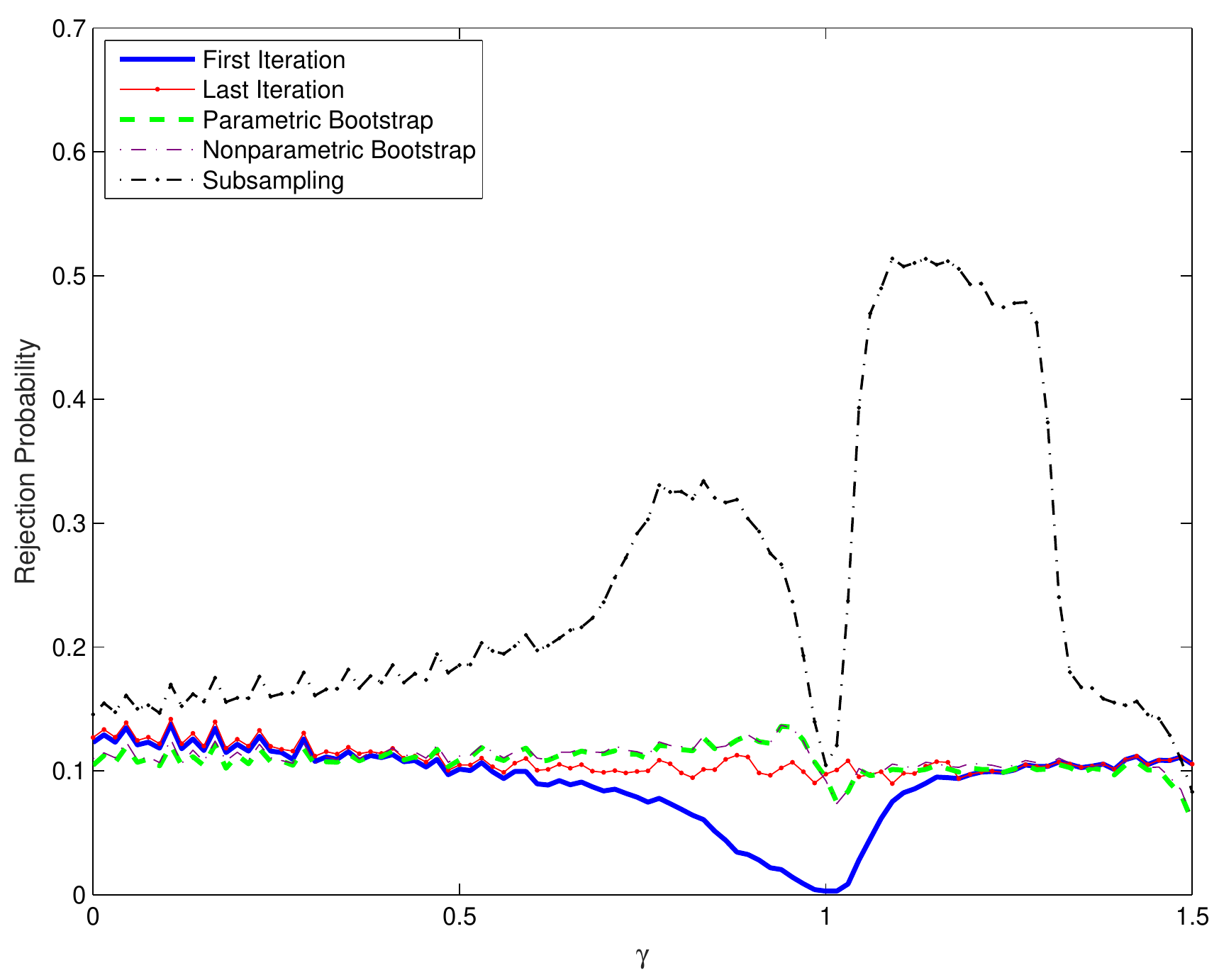}
\subcaption{$\rho = 0.95$} \endminipage\hfill \minipage{0.5\textwidth} %
\centering %
\includegraphics[width=6.8cm]{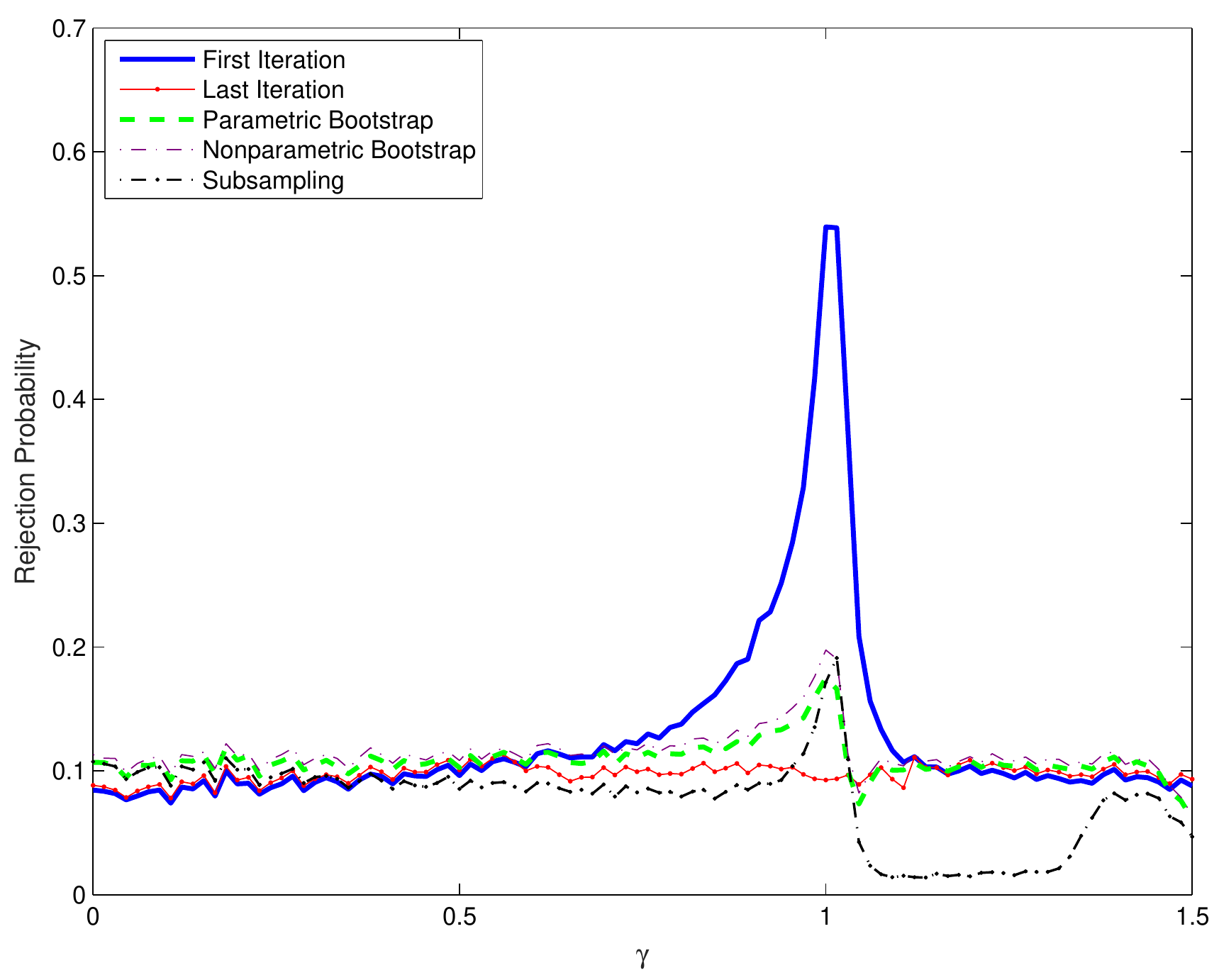}
\subcaption{$\rho = -0.95$} \endminipage\hfill
\end{figure}

Figure \ref{fig:Size} plots null rejection probabilities for $\rho =0.95$
and $\rho =-0.95$. Taking $c$ ranging from $-100$ to $50$ allows us to
consider cases in which the regressors are stationary ($\gamma
=0.00,...,0.50 $), nearly integrated ($\gamma =0.85,...,0.95$), exactly
integrated ($\gamma =1.00$), nearly explosive ($\gamma =1.05$ or $1.10$),
and explosive ($\gamma =1.20,...,1.50$).

When $\rho =0.95$ and we include only the initial points $\gamma =0.5$ and $%
\gamma =1.2$, the null rejection probability is lower than 10\% for values
of $\gamma $ between the endpoints $\gamma =0.5$ and $\gamma =1.2$. The
associated test seems to have correct size but the null rejection
probability is near zero when $\gamma =1$. As a result, the power of a test
using only the two endpoints can be very low when the series is nearly
integrated. By adding additional points, the null rejection probability
becomes closer to 10\% for all values $\gamma \in \left[ 0.5,1.2\right] $.
In practice, our algorithm includes about 10 additional points. When $\rho
=-0.95$ and we consider only the initial points $\gamma =0.5$ and $\gamma
=1.2$, the null rejection probabilities can be as large as 55\%. Applied
researchers would then reject the null considerably more often than the
reported nominal level. By sequentially adding the most discrepant points,
the algorithm eventually yields a null rejection probability close to 10\%.

Figure \ref{fig:Size} also compares the CVF approach with the usual
bootstrap method, the parametric bootstrap, and subsampling. For the
(nonparametric) bootstrap, we draw a random sample from the centered
empirical distribution function of the fitted residuals for $\epsilon _{t}$
from OLS regressions. For the parametric bootstrap, we replace the unknown
autoregressive parameter $\gamma $ with the OLS\ estimator. For the
subsampling, we choose the subsampling parameters following\ the
recommendation of \citet{RomanoWolf01}. The bootstrap sampling schemes have
similar performance and do not have null rejection probability near 10\%
when the series are nearly integrated. This failure is due to the fact that
the estimator $\widehat{\beta }$ is not locally asymptotically equivariant
in the sense of \citet{Beran97}; see also \citet{BMMRT91}. The subsampling
is a remarkably general method and has good performance when $\gamma =1$.
However, its null rejection probability can be quite different from 10\% for
values of $\gamma $ near one. In the supplement, we show that the bootstrap
and subsampling can be even further away from the nominal level of 10\% if
we make inference in the presence of a time trend.

\bigskip

\begin{figure}[tbh]
\caption{Critical Value Function}
\label{fig:CVF}
\bigskip \minipage{0.5\textwidth} \centering %
\includegraphics[width=6.8cm]{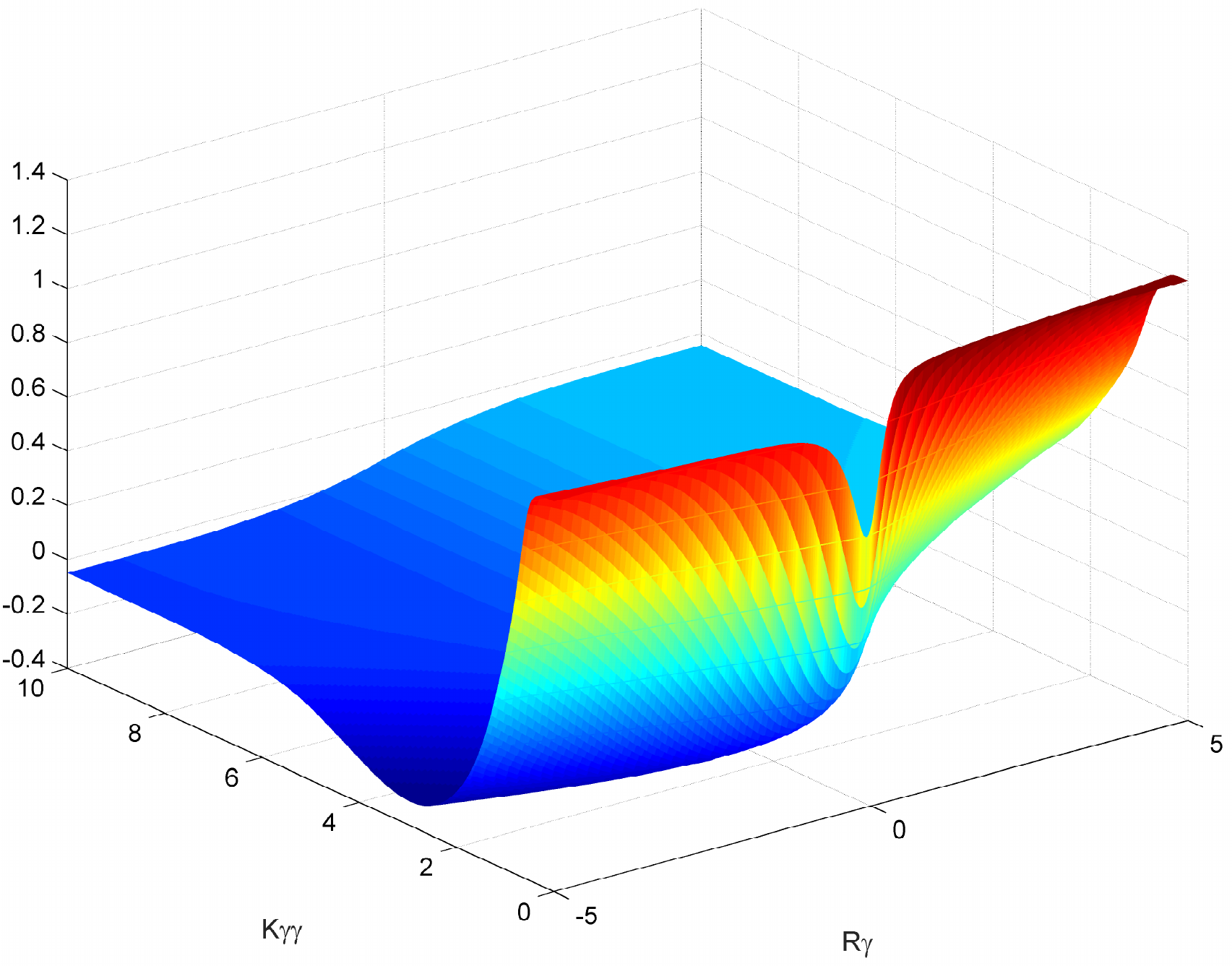}
\subcaption{$\rho = 0.95$} \endminipage \hfill \minipage{0.5\textwidth} %
\centering %
\includegraphics[width=6.8cm]{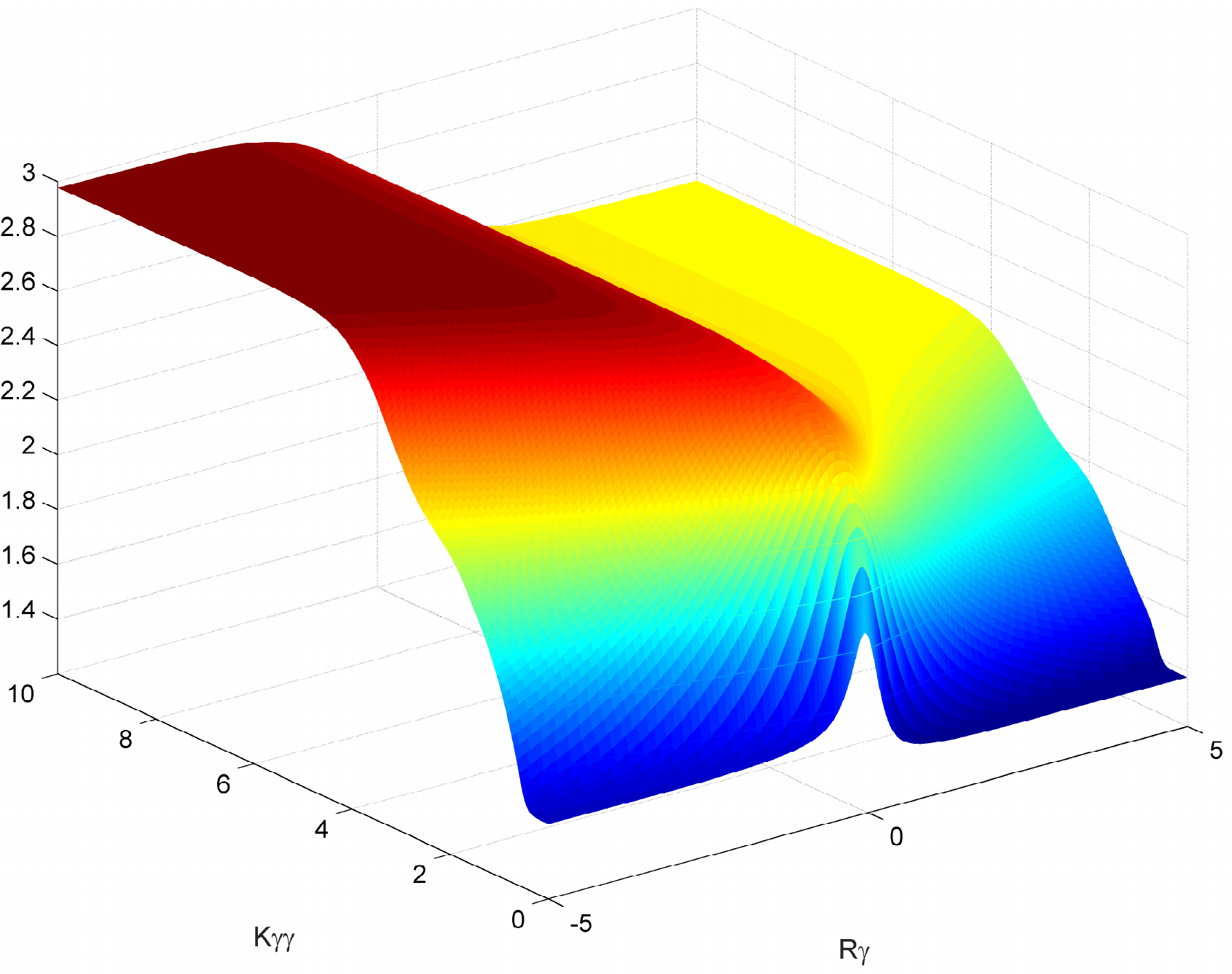}
\subcaption{$\rho = -0.95$} \endminipage \hfill
\end{figure}

Interestingly, the null rejection probability for the CVF method is not far
from 10\% even for values of $\gamma $ outside the range $\left[ 0.5,1.2%
\right] $ whether $\rho $ is positive or negative. This suggests that the
critical value functions take values close to 1.28, the 90\% quantile of a
standard normal distribution, when the series is stationary or explosive.
Figure \ref{fig:CVF} supports this observation when the series is
stationary. This graph plots the critical value function as a mapping of $%
R_{\gamma ,T}\left( 1\right) $ and $K_{\gamma \gamma ,T}\left( 1\right) $.
To shorten notation, we denote these statistics $R_{\gamma }$ and $K_{\gamma
\gamma }$ respectively. Standard asymptotic arguments show that $K_{\gamma
\gamma }$ converges in probability to zero and $R_{\gamma }/K_{\gamma \gamma
}$ diverges to $-\infty $. Whether $\rho =0.95$ or $\rho =-0.95$, the
critical value function takes values close to 1.28 when $K_{\gamma \gamma }$
is close to zero and $R_{\gamma }$ is bounded away from zero (which, of
course, implies that their ratio diverges). We cannot visualize, in this
figure, the explosive case in which $K_{\gamma \gamma }$ diverges to $\infty
$. In the supplement, we plot the critical values as a function of $%
R_{\gamma }/K_{\gamma \gamma }$ and $K_{\gamma \gamma }$ and re-scale the
axis to accommodate for the explosive case as well. We show again that the
critical value function takes values near 1.28 also when the series is
explosive.

\subsection{Size and Power Comparison\label{Power SubSec}}

We now evaluate different size correction methods for the t-statistic when
the variance matrix for the innovations $\varepsilon _{t}^{y}$ and $%
\varepsilon _{t}^{x}$ is unknown. The covariance matrix is estimated using
residuals from OLS regressions.

We compare size and power of the similar t-test with the two other similar
tests: the $L_{2}$ test of \citet{Wright00} and the UMPCU\ test of %
\citet{JanssonMoreira06}. The figures present power curves for $c=-15$ and $%
c=0$. The rejection probabilities are plotted against local departures from
the null hypothesis. Specifically, the parameter of interest is $\beta
=b\cdot \sigma _{yy.x}\cdot g_{T}\left( \gamma \right) $ for $%
b=-10,-9,...,10 $. The value $b=0$ corresponds to $\beta =0$ for the null
hypothesis $H_{0}:\beta \leq 0$.

\begin{figure}[tbh]
\caption{Power ($\protect\rho =0.95$)}
\label{fig:Power positive}
\bigskip \centering \minipage{0.5\textwidth} %
\includegraphics[width=6.8cm]{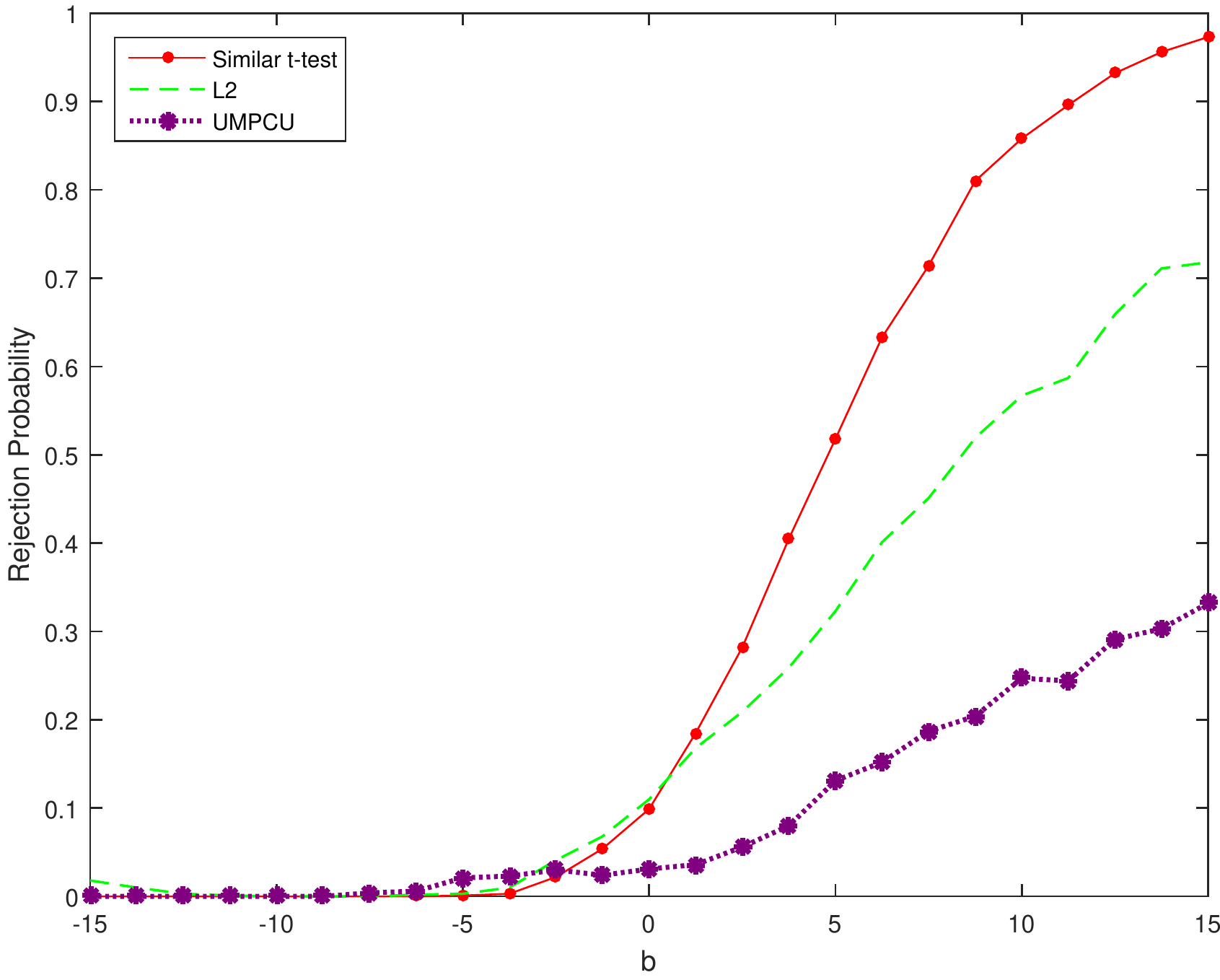}
\subcaption{$\gamma = 0.85$} \endminipage\hfill \minipage{0.5\textwidth} %
\includegraphics[width=6.8cm]{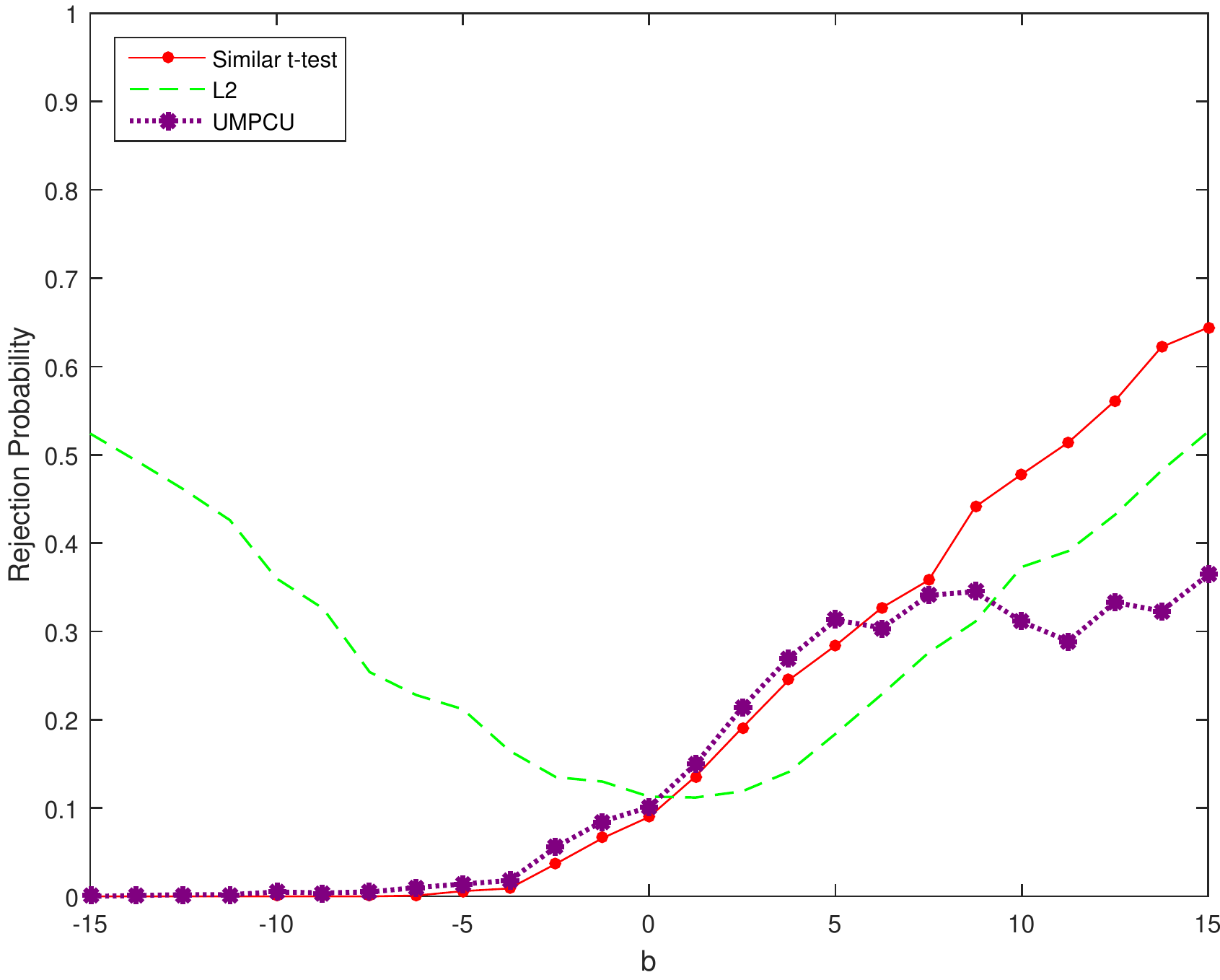}
\subcaption{$\gamma = 1$} \endminipage \hfill
\end{figure}

Figure \ref{fig:Power positive} plots size and power comparisons for $\rho
=0.95$. The similar t-test has correct size and null rejection probability
near the nominal level when $b=0$. Hence, replacing the unknown variance by
a consistent estimate has little effect on size in practice. Indeed, the
limiting distribution of the t-statistic and the CVF is the same whether the
variance is estimated or not\footnote{%
More generally, it is even possible to accommodate for heteroskedasticity
and autocorrelation for the innovations by adjusting the statistics $R_{T}$
and $K_{T}$ and estimating a long-run variance, as done by %
\citet{JanssonMoreira06}.}. The $L_{2}$ test is also similar but does not
have correct size when the series is integrated. The $L_{2}$ test actually
has a bell-shaped power curve when the regressor is integrated, but has a
monotonic power curve when the series is nearly integrated. This
inconsistency makes it hard to use the $L_{2}$ test for either one-sided or
two-sided hypothesis testing problems. Although the UMPCU test is a similar
test, our reported null rejection probability is near zero when $b=0$ and
the series is nearly integrated. This problem is possibly due to the
inaccuracy of the algorithm used by \citet{JanssonMoreira06} for integral
approximations (and their computed critical value).

Whether the series is nearly integrated or exactly integrated, the similar
t-test dominates the two other similar tests in terms of power. When $c=-15$
and $b=10$, the similar t-test has power above 90\% while the $L_{2}$ and
UMPCU tests have power lower than 60\% and 30\%, respectively. When $c=0$
and $b=10$, the similar t-test has power close to 50\% while the $L_{2}$ and
UMPCU\ tests have power close to 40\% and 30\%. The fact that the similar
t-test outperforms the UMPCU test may seem puzzling at first. However, the
UMPCU\ test is optimal within the class of similar tests conditional on
specific ancillary statistics $S_{\beta \beta }$ and $S_{\gamma \gamma }$.
As there exist similar tests beyond those considered by %
\citet{JanssonMoreira06}, our less restrictive requirement yields power
gains. Therefore, the performance of the similar t-test shows that the
UMPCU\ test uses a conditional argument that entails unnecessary power loss.

The similar t-test also has considerably better performance than the other
similar tests when the endogeneity coefficient is negative. Figure \ref%
{fig:Power negative} presents power curves for $\rho =-0.95$. All three
tests present null rejection probabilities near the nominal level when $b=0$%
. The $L_{2}$ test again can present bell-shaped power, and does not even
have correct size. The power functions for the similar t-test and UMPCU\
test are monotonic and have null rejection probabilities smaller than 10\%
when $c=-15$ or $c=0$.
\begin{figure}[tbh]
\caption{Power ($\protect\rho =-0.95$)}
\label{fig:Power negative}
\bigskip \centering \minipage{0.5\textwidth} %
\includegraphics[width=6.8cm]{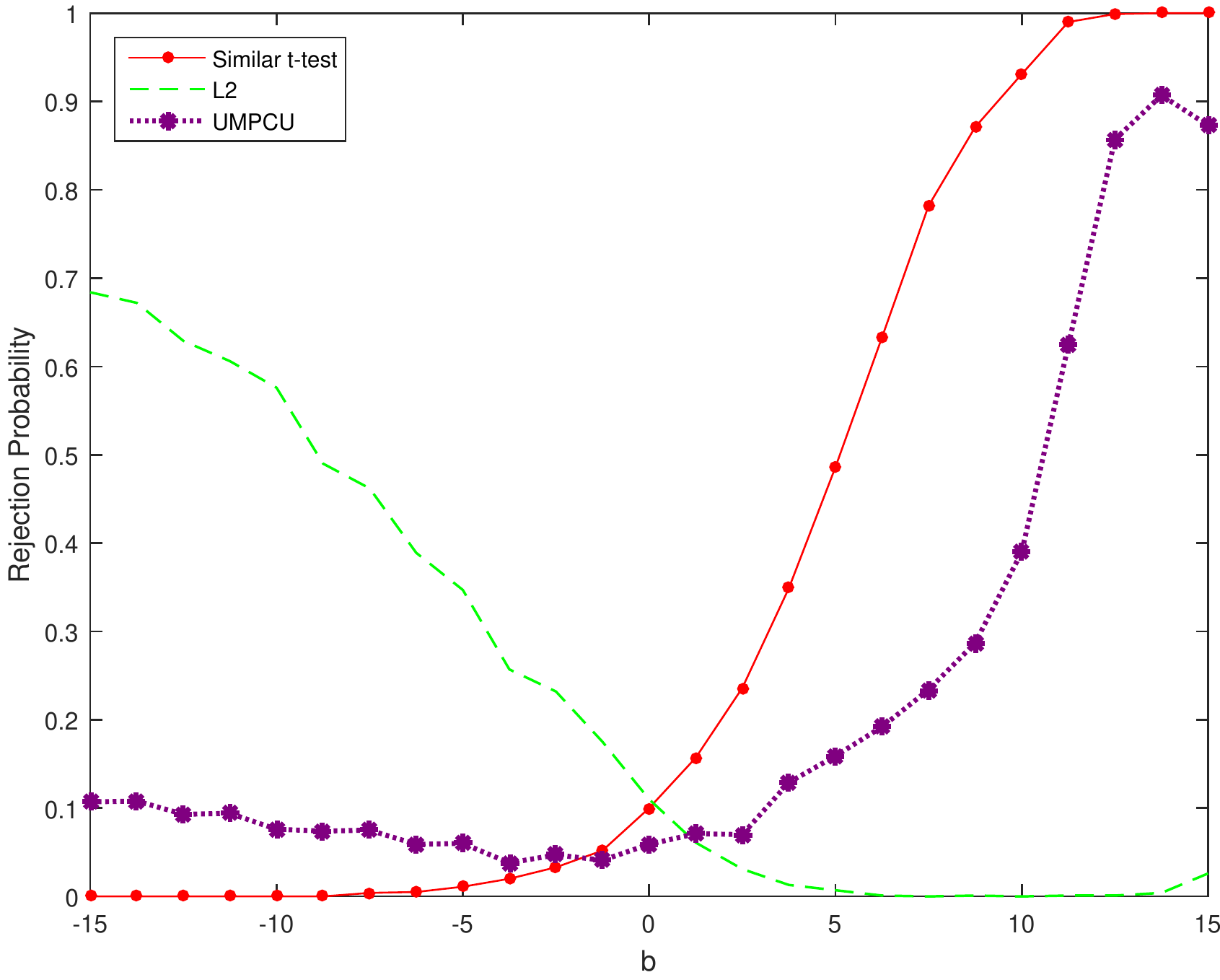}
\subcaption{$\gamma = 0.85$} \endminipage\hfill \minipage{0.5\textwidth} %
\includegraphics[width=6.8cm]{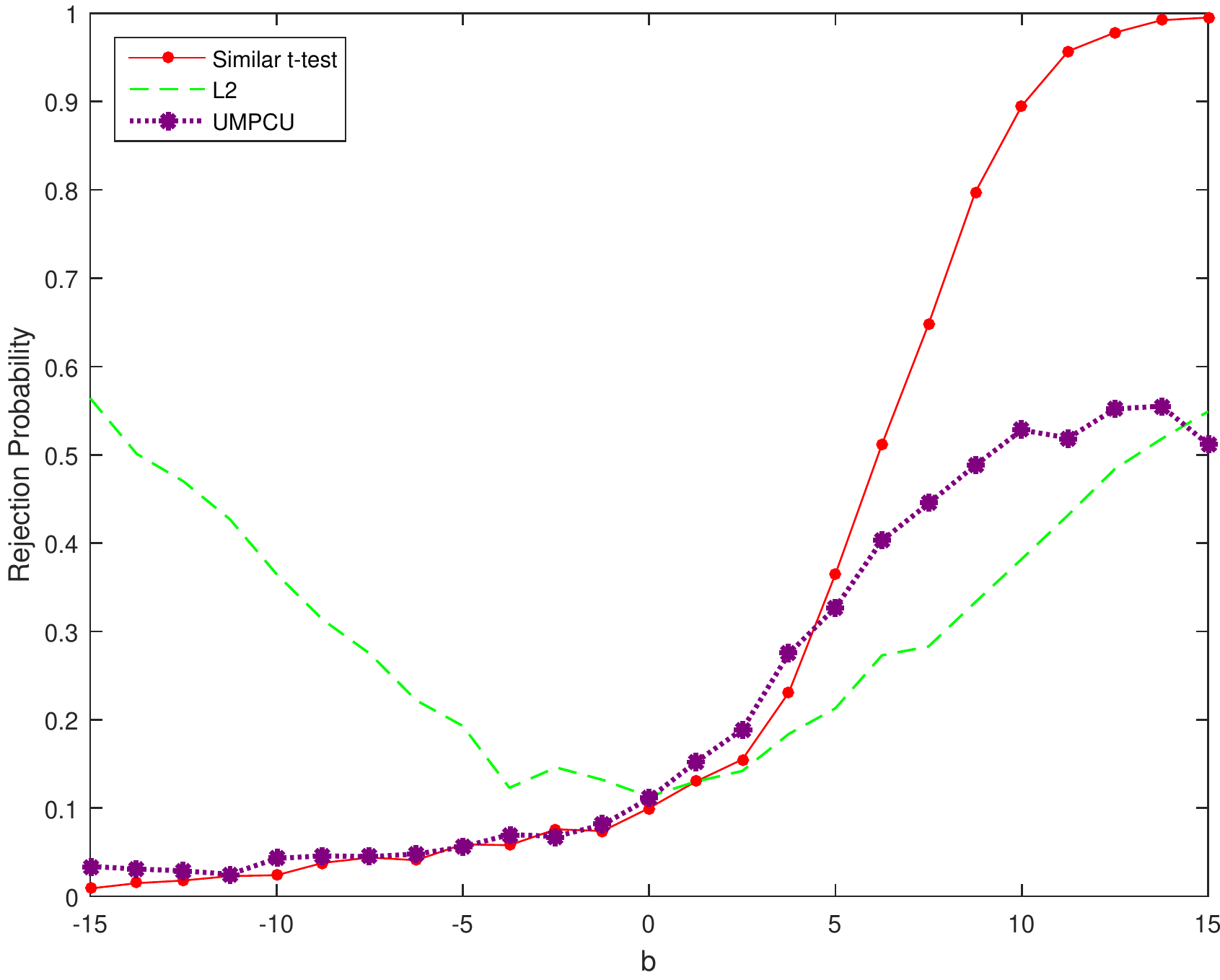}
\subcaption{$\gamma = 1$} \endminipage \hfill
\end{figure}

The supplement presents size and power comparisons for all combinations of $%
c=-100$, $-50$, $-15$, $0$, $10$, $30$ and $\rho =0.5$, $-0.5$, $0.95$, $%
-0.95$. These numerical simulations further show that the similar t-test has
overall better performance than the $L_{2}$ and UMPCU\ tests in terms of
size and power.

\section{Conclusion and Extensions\label{Conclusion Sec}}

This paper proposes a novel method to find tests with correct size for
situations in which conventional numerical methods do not perform
satisfactorily. The critical value function (CVF) approach is very general
and relies on weighting schemes which can be found through linear
programming (LP). Considering the polynomial speed of interior-point methods
for LP, it is fast and straightforward to find the CVF. The weights are
chosen so that the test is similar for a finite number of null rejection
probabilities. If the nuisance parameter dimension is not too large, we
expect the CVF method to work adequately.

In a model with persistent regressors, the CVF method yields a similar
t-test which outperforms other similar tests proposed in the literature. It
would be interesting to assess how our methodology performs in other models.
For example, take the autoregressive model. If the model has one lag, %
\citet{Hansen99} suggests a \textquotedblleft grid\textquotedblright\
bootstrap which delivers confidence regions with correct coverage
probability (in a uniform sense). This method relies on less restrictive
conditions than the bootstrap. However, it requires the null distribution of
the test statistic to be asymptotically independent of the nuisance
parameter (hence, it is applicable even for our model with persistent
regressors). If the autoregressive model has more lags, \citet{RomanoWolf01}
propose subsampling methods based on Dickey-Fuller representations. Their
proposed confidence region which has correct coverage level can be
conservative (in the sense that some of the coverage probabilities can be
larger than the nominal level). By inverting the tests based on the CVF
method, we are able to provide confidence regions which instead have
coverage probabilities near the nominal level.

\bibliographystyle{Chicago}
\bibliography{References}



\end{document}


\setcounter{table}{0} \renewcommand{\thetable}{S\arabic{table}} %
\setcounter{figure}{0} \renewcommand{\thefigure}{S\arabic{figure}}

\begin{frontmatter}

\title{Supplement to ``A Critical Value Function Approach, with an
Application to Persistent Time-Series''}
\runtitle{A Critical Value Function Approach}


\begin{aug}
\author{\fnms{Marcelo J.} \snm{Moreira}\ead[label=e1]{moreiramj@fgv.br}},
and
\author{\fnms{Rafael} \snm{Mour\~{a}o}\ead[label=e2]{rmrodrigues@fgvmail.br}},

\runauthor{Moreira, Mour\~{a}o, and Moreira}

\affiliation{FGV/EPGE}

\address{Escola de P\'{o}s-Gradua\c{c}\~{a}o em Economia e Finan\c{c}as (FGV/EPGE)\\
Getulio Vargas Foundation - 11th floor\\
Praia de Botafogo 190\\
Rio de Janeiro - RJ 22250-040\\
\printead{e1}}

\end{aug}

\end{frontmatter}


\section{Introduction}

This supplement consists of two parts. The first section includes additional
graphs on the critical value function, size and power comparison. The second
section provides additional proofs for the remaining theoretical results.

\section{Numerical Results}

We provide additional numerical results. First, we allow for a time trend
and compare the performance of commonly-used sampling methods to our CVF
approach. We then consider flattening schemes to guarantee that the critical
value function (CVF) converges to the usual standard normal quantiles.
Finally, we make further power comparisons between our similar t-test and
other similar tests for different values of the autoregressive parameter $%
\gamma $ and degrees of endogeneity $\rho $.

\subsection{Time Trend}

For now, we extend the model to allow for a time trend. The model is given by%
\begin{eqnarray*}
y_{t} &=&\mu _{y}^{\prime }d_{t}+\beta x_{t-1}+\epsilon _{t}^{y} \\
x_{t} &=&\gamma x_{t-1}+\epsilon _{t}^{x},
\end{eqnarray*}%
where $\mu _{y}$ is a two-dimensional vector and $d_{t}=\left( 1,t\right)
^{\prime }$. We consider the group of translation transformations on the data%
\begin{equation*}
g\circ \left( y,x\right) =\left( y+Dg,x\right) ,
\end{equation*}%
where $g$ is a two-dimensional vector and $D$ is a $T\times 2$ matrix whose $%
t$-th line is given by $d_{t}^{\prime }$. This yields a transformation on
the parameter space%
\begin{equation*}
g\circ \left( \beta ,\gamma ,\mu _{y}\right) =\left( \beta ,\gamma ,\mu
_{y}+g\right) .
\end{equation*}

Any invariant test can be written as a function of the maximal invariant
statistic. Let $q=\left( q_{1},q_{2}\right) $ be an orthogonal $T\times T$
matrix where the first column is given by $q_{1}=D\left( D^{\prime }D\right)
^{-1/2}$. Algebraic manipulations show that $q_{2}q_{2}^{\prime }=M_{D}$,
where $M_{D}=I_{T}-D\left( D^{\prime }D\right) ^{-1}D^{\prime }$ is the
projection matrix to the space orthogonal to $1_{T}$. Let $x_{-1}$ be the $T$%
-dimensional vector whose $t$-th entry is $x_{t-1}$, and define $w^{\mu
}=q_{2}^{\prime }w$ for a $T$-dimensional vector $w$. The maximal invariant
statistic is given by $r=\left( y^{\mu },x\right) $. Its density function is
given by
\begin{eqnarray}
&&f_{\beta ,\gamma }\left( y^{\mu },x\right) =\left( 2\pi \sigma
_{xx}\right) ^{-\frac{T}{2}}\exp \left\{ -\frac{1}{2\sigma _{xx}}%
\sum_{t=1}^{T}\left( x_{t}-x_{t-1}\gamma \right) ^{2}\right\}
\label{(likelihood eq)} \\
&&\times \left( 2\pi \sigma _{yy.x}\right) ^{-\frac{T-2}{2}}\exp \left\{ -%
\frac{1}{2\sigma _{yy.x}}\sum_{t=1}^{T}\left( y_{t}^{\mu }-x_{t}^{\mu }\frac{%
\sigma _{xy}}{\sigma _{xx}}-x_{t-1}^{\mu }\left[ \beta -\gamma \frac{\sigma
_{xy}}{\sigma _{xx}}\right] \right) ^{2}\right\} ,  \notag
\end{eqnarray}%
where $\sigma _{yy.x}=\sigma _{yy}-\sigma _{xy}^{2}/\sigma _{xx}$ is the
variance of $\epsilon _{t}^{y}$ not explained by $\epsilon _{t}^{x}$.
Because of invariance, the density depends on $\beta $ and $\gamma $, but
not on $\mu _{y}$. Furthermore, the formulae for the sufficient statistics
and the log-likelihood ratio process are exactly the same as in the model
with only an intercept, but replacing $q_{2}$ in $y^{\mu }$ and $x^{\mu }$.
The asymptotic distribution of the new $R_{T}\left( \gamma \right) $ and $%
K_{T}\left( \gamma \right) $ statistics is, however, different when we allow
for time trend. In particular, the model curvature increases in the sense of %
\citet{Efron75,Efron78} when $\gamma =1$. This yields larger departures from
the nominal size by conventional sampling schemes.

\begin{figure}[tbh]
\caption{Size (intercept and time trend)}
\label{fig:Size trend}
\bigskip \minipage{0.5\textwidth} \centering %
\includegraphics[width=6.8cm]{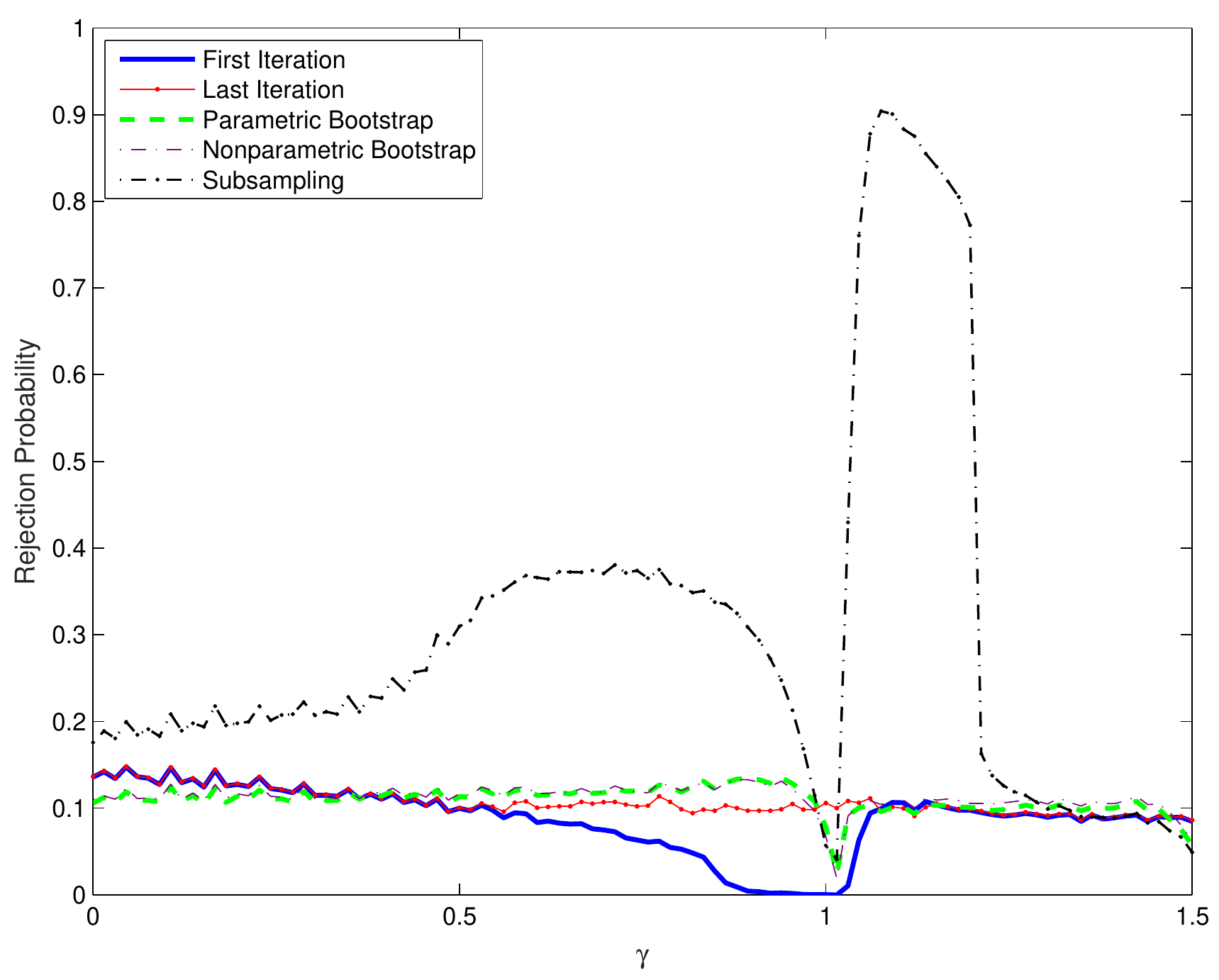}
\subcaption{$\rho = 0.95$} \endminipage\hfill \minipage{0.5\textwidth} %
\centering %
\includegraphics[width=6.8cm]{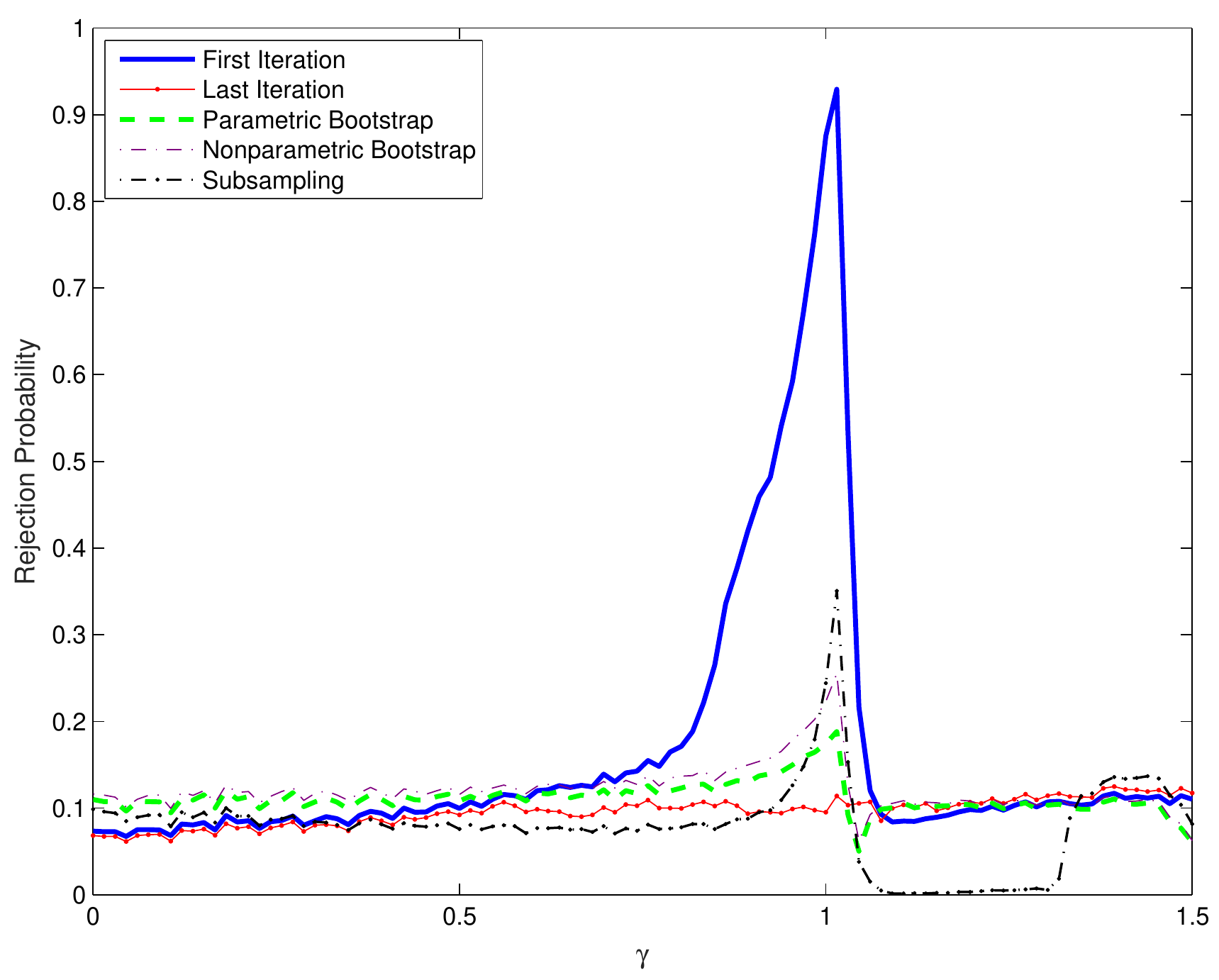}
\subcaption{$\rho = -0.95$} \endminipage\hfill
\end{figure}
Figure \ref{fig:Size trend} plots null rejection probabilities when we allow
for a time trend. Comparing this graph with Figure 1, it is
evident that both bootstrap methods and subsampling have null rejection
probabilities farther away from the 10\% nominal size when we allow for a
time trend. The CVF, however, is able to correct size distortions whether
there is a time trend or not.
\begin{figure}[tbh]
\caption{Critical Value Function (intercept and time trend)}
\label{fig:CVF trend}
\bigskip \minipage{0.5\textwidth} \centering %
\includegraphics[width=6.8cm]{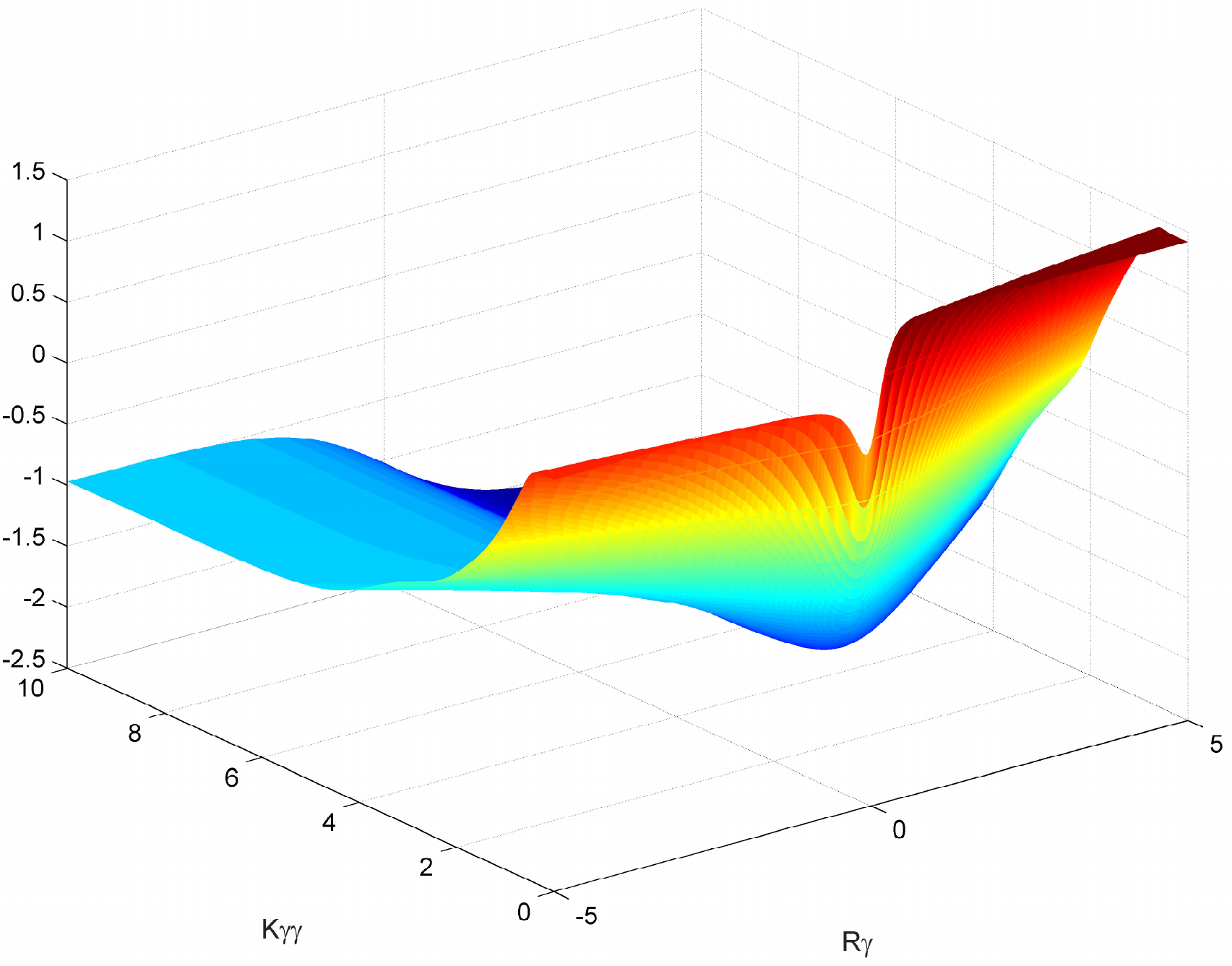}
\subcaption{$\rho = 0.95$} \endminipage \hfill \minipage{0.5\textwidth} %
\centering %
\includegraphics[width=6.8cm]{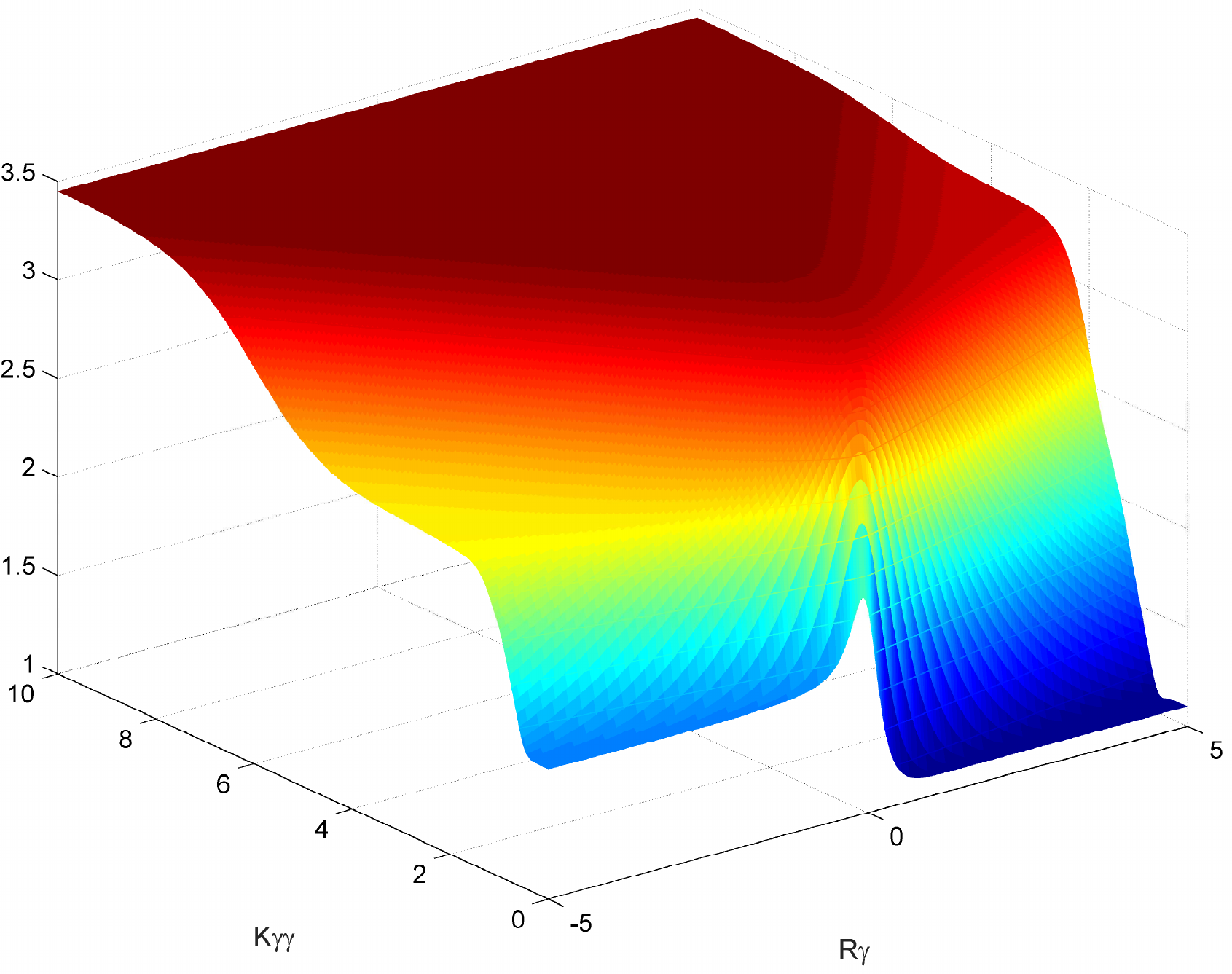}
\subcaption{$\rho = -0.95$} \endminipage \hfill
\end{figure}

As such, the CVF needs to provide more curvature when we allow for time
trend. This can be seen by comparing Figure \ref{fig:CVF trend}, which plots
the CVF when the model has a time trend, to Figure 2. The range
of values that the CVF can take is larger when we make inference in the
presence of a time trend. For example, take the CVF plot when $\rho =0.95$.
The CVF can have values as low as -2.5 and as high as 1.5 in the presence of
a time trend. On the other hand, the CVF ranges from -1 to 1.5 when we make
inference knowing there is no time trend.

\subsection{Flattening Schemes}

We study the behavior of the critical value function (CVF) when the series
is stationary and explosive. We focus on the main model in which there is
only intercept and no time trend. We are able to analyze the behavior of $%
R_{\gamma }/K_{\gamma \gamma }$ and $K_{\gamma \gamma }$ when the series is
stationary and explosive (to save space, we again use the notation $%
R_{\gamma }$ and $K_{\gamma \gamma }$ for $R_{\gamma ,T}\left( 1\right) $
and $K_{\gamma \gamma ,T}\left( 1\right) $, respectively). Figure \ref%
{fig:CVF suppl} plots the critical values as a function of $R_{\gamma
}/K_{\gamma \gamma }$ and $K_{\gamma \gamma }$. It is convenient to re-scale
the x-axis and y-axis in the figures to analyze the critical value under
different asymptotic regimes. We choose the scale
\begin{equation*}
G\left( z\right) \equiv 2.\left( F\left( z\right) -0.5\right) \text{,}
\end{equation*}%
where $F\left( z\right) =\left[ 1+\exp \left( -z\right) \right] ^{-1}$ is
the logistic distribution. This scale is an increasing, symmetric
transformation, which preserves zero and maps $\infty $ into $1$; that is, $%
G^{\prime }\left( z\right) >0$,~$G\left( -z\right) =-G\left( z\right) $, $%
G\left( 0\right) =0$, and lim$_{z\rightarrow \infty }G\left( z\right) =1$.
\begin{figure}[tbh]
\caption{Critical Value Function}
\label{fig:CVF suppl}
\bigskip \minipage{0.5\textwidth} \centering %
\includegraphics[width=6.8cm]{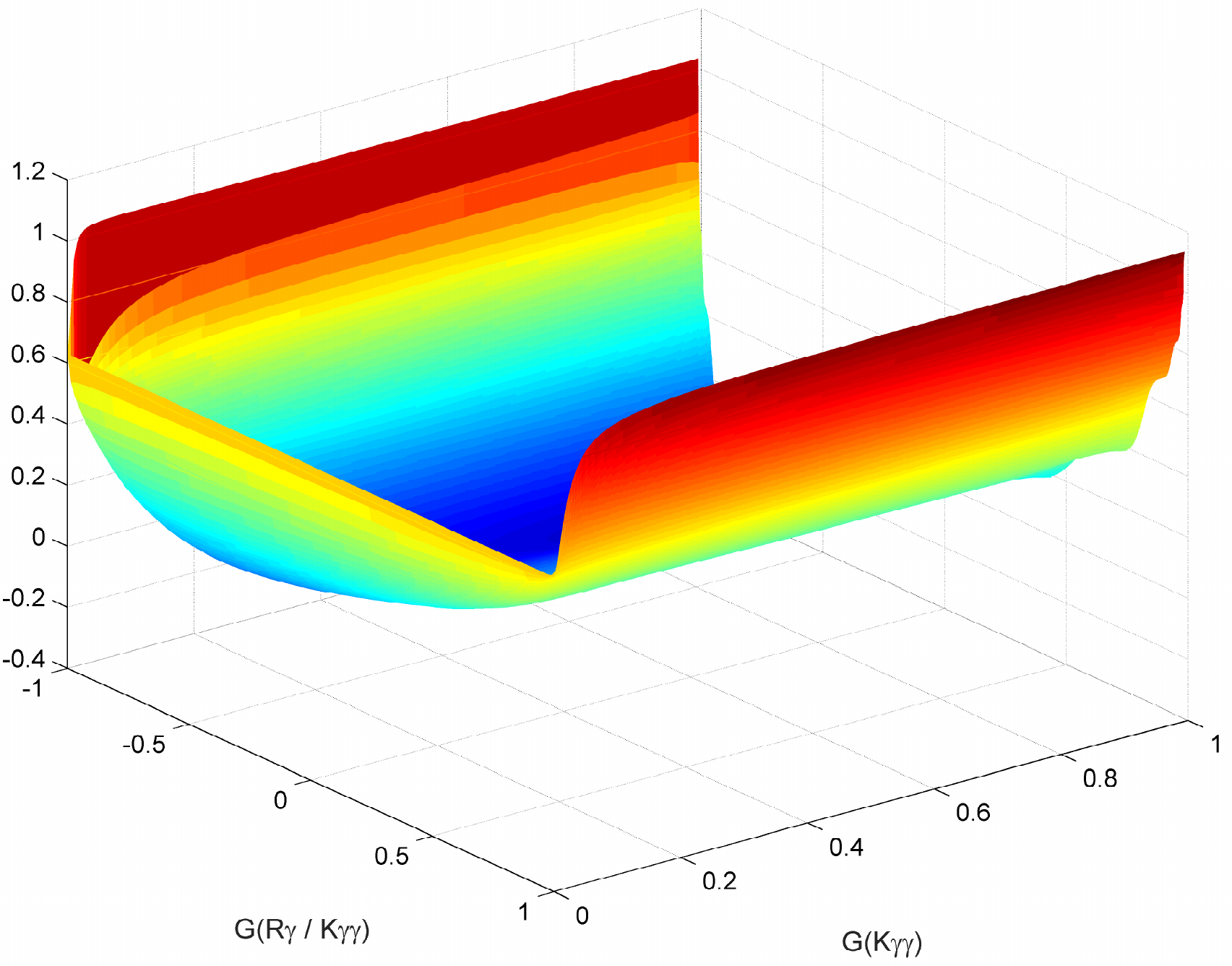}
\subcaption{$\rho = 0.95$} \endminipage \hfill \minipage{0.5\textwidth} %
\centering %
\includegraphics[width=6.8cm]{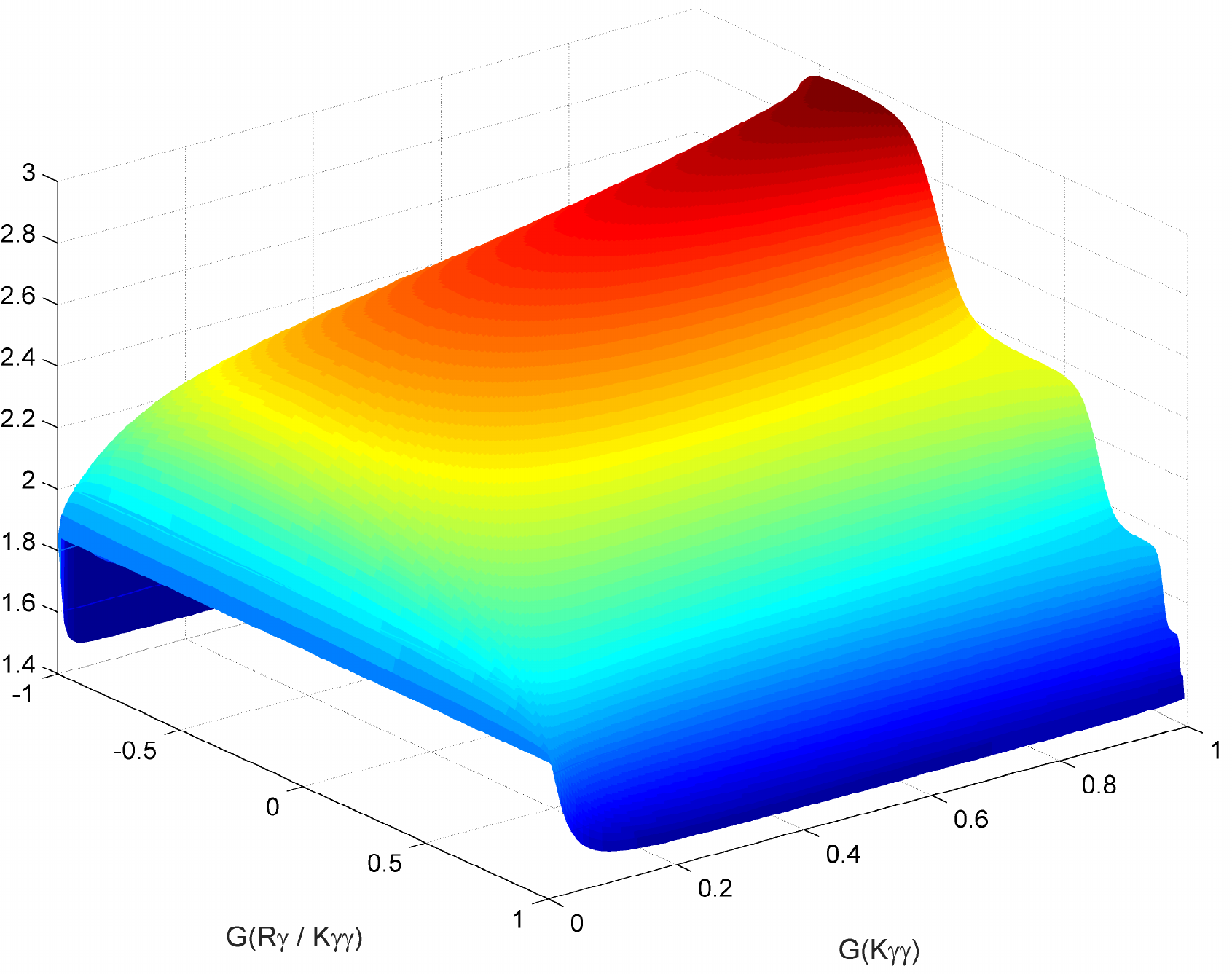}
\subcaption{$\rho = -0.95$} \endminipage \hfill
\end{figure}

\begin{figure}[tbh]
\caption{Size}
\label{fig:Size suppl}
\bigskip \minipage{0.5\textwidth} \centering %
\includegraphics[width=6.8cm]{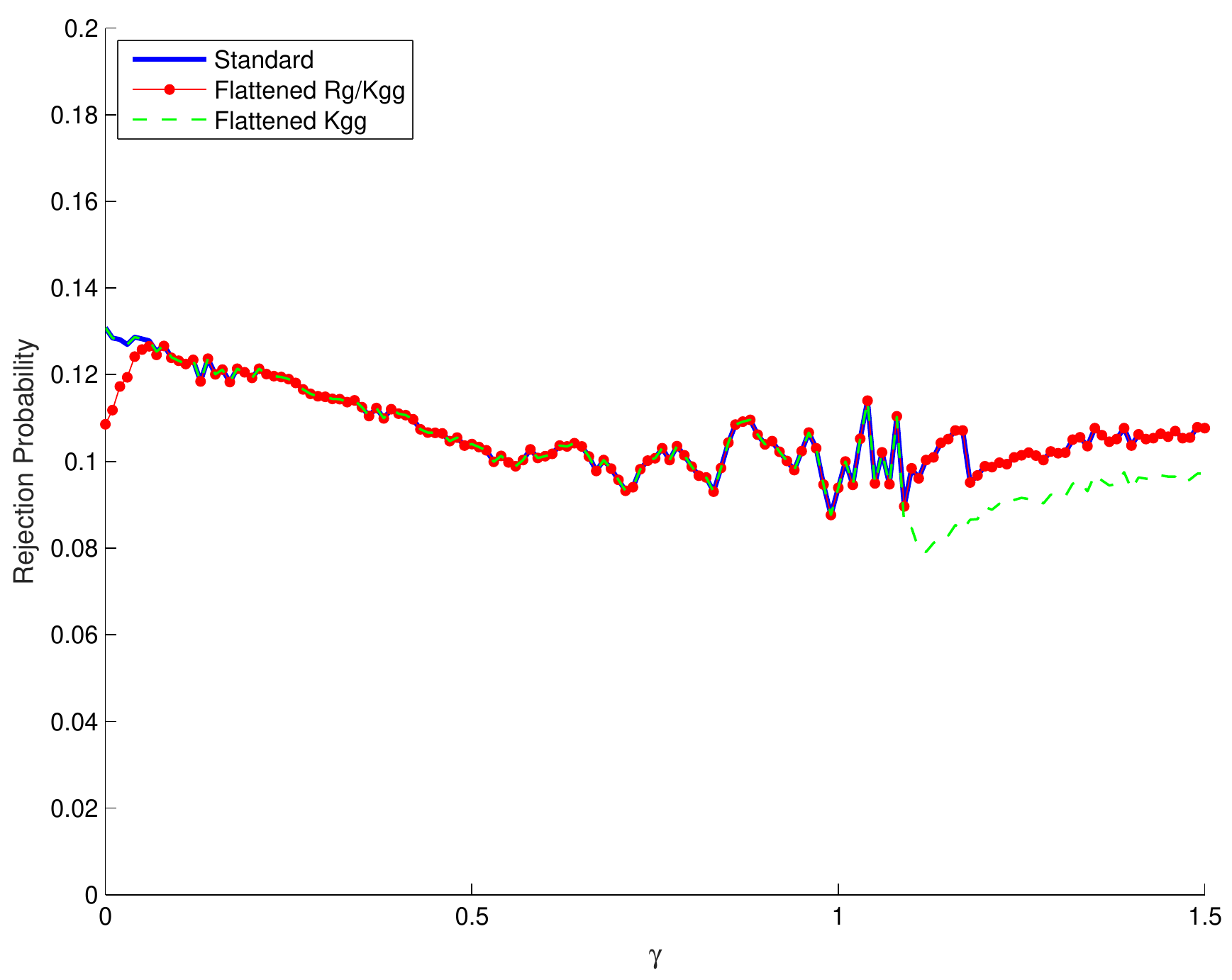}
\subcaption{$\rho = 0.95$} \endminipage\hfill \minipage{0.5\textwidth} %
\centering %
\includegraphics[width=6.8cm]{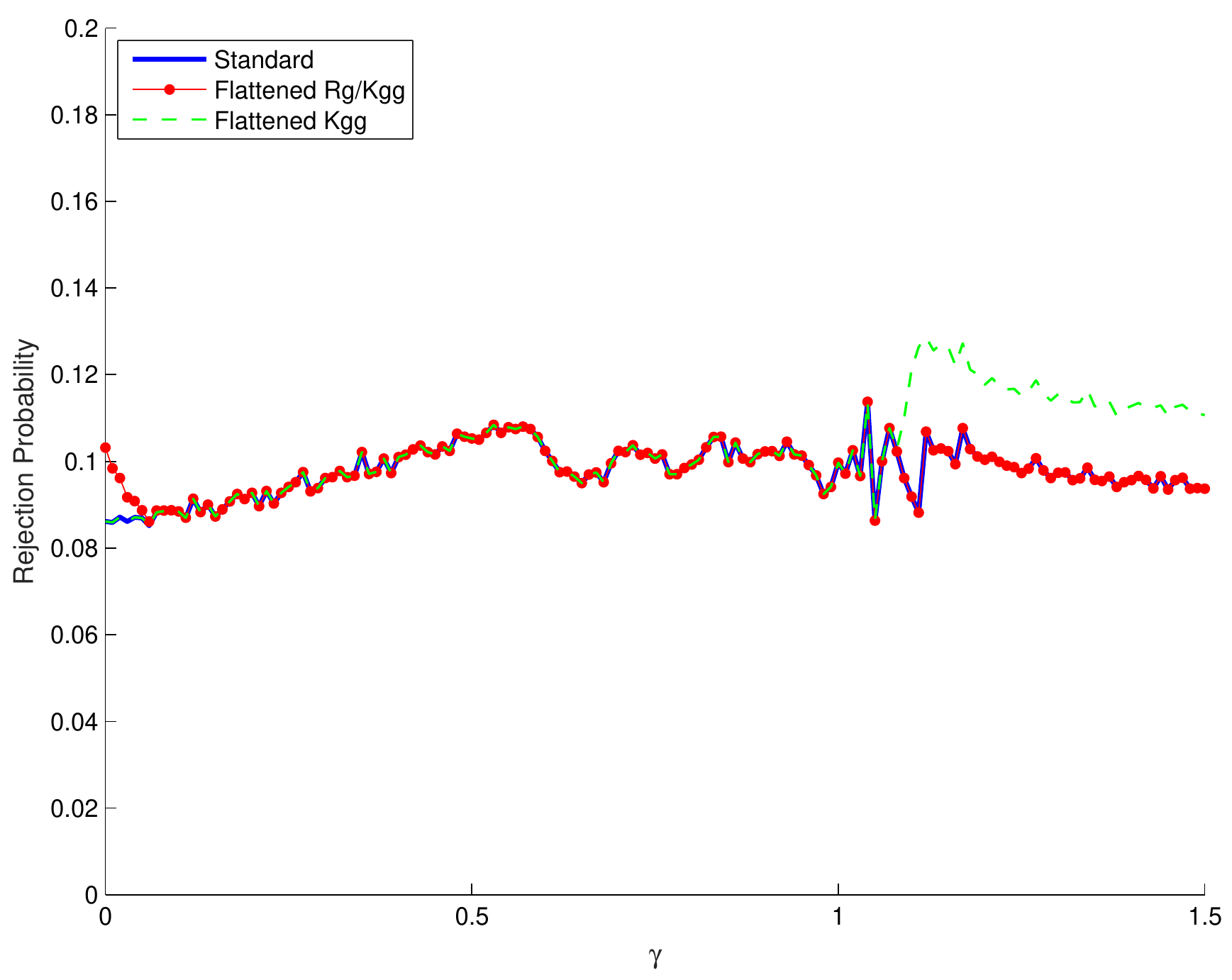}
\subcaption{$\rho = -0.95$} \endminipage\hfill
\end{figure}

When the series is stationary, $R_{\gamma }/K_{\gamma \gamma }$ converges to
$-\infty $ (under local alternatives) and $K_{\gamma \gamma }$ to $0$. This
convergence is represented by the point $\left( -1,0\right) $ in the domain
for $G\left( R_{\gamma }/K_{\gamma \gamma }\right) $ and $G\left( K_{\gamma
\gamma }\right) $. When the series is explosive, $R_{\gamma }/K_{\gamma
\gamma }$ diverges to $\infty $ (under local alternatives) and $K_{\gamma
\gamma }$ to $\infty $. This limit is indicated by the point $\left(
1,1\right) $ in the new domain. In both situations, the critical value
function takes values near 1.28, the 90\% quantile of a standard normal
distribution. Of course, those values are not exactly equal to 1.28. This
suggests we can try to flatten the critical value function using the values
of $R_{\gamma }/K_{\gamma \gamma }$ or $K_{\gamma \gamma }$ themselves. For
example, we could consider setting the critical value function at 1.28 for
(i) large values of $\left\vert R_{\gamma }/K_{\gamma \gamma }\right\vert $,
or (ii) small values and large values of $K_{\gamma \gamma }$. The figures
below consider flattening schemes when (i) $\left\vert R_{\gamma }/K_{\gamma
\gamma }\right\vert >10^{2}$ and (ii) $K_{\gamma \gamma }<10^{-2}$ or $%
K_{\gamma \gamma }>10^{6}$. For these regions, the critical value function
is forced to take the value 1.28. Figure \ref{fig:Size suppl} plots the null
rejection probabilities from both flattening methods for $\rho =0.95$ and $%
-0.95$.

The tests' sizes are sometimes sensitive to the choice of the thresholds for
both flattening methods. When flattening slightly improves size, it does not
yield a substantial effect on the power of the respective tests. Hence, we
report power using the standard critical value function, without any
auxiliary flattening method.

\subsection{Additional Power Comparison}

\begin{figure}[tbh]
\caption{Power ($\protect\rho =0.5$)}
\label{fig:Power 0.5 suppl}
\bigskip \minipage{0.5\textwidth} \centering %
\includegraphics[width=6.8cm]{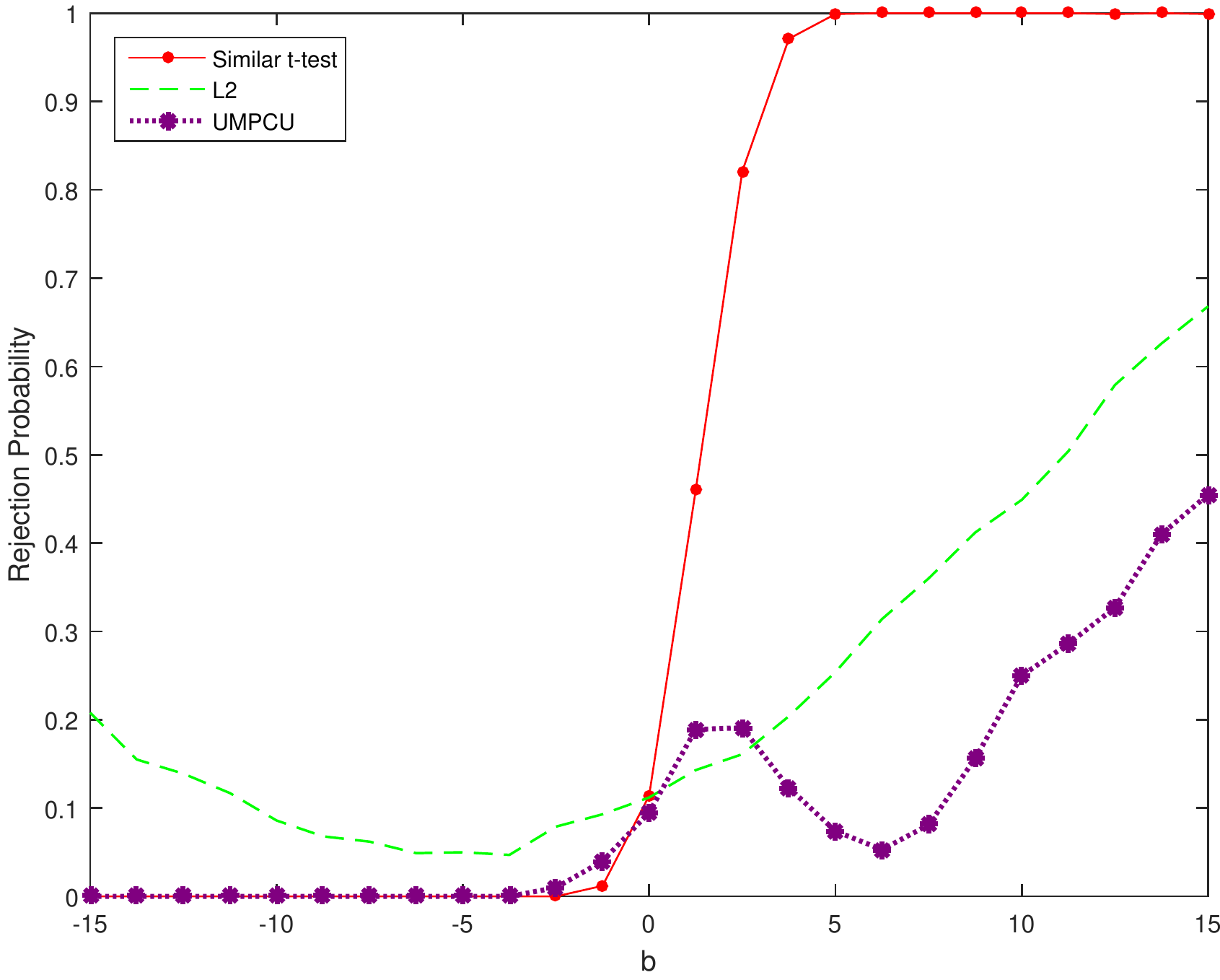} %
\subcaption{$\gamma = 0$} \endminipage \hfill \minipage{0.5\textwidth} %
\centering %
\includegraphics[width=6.8cm]{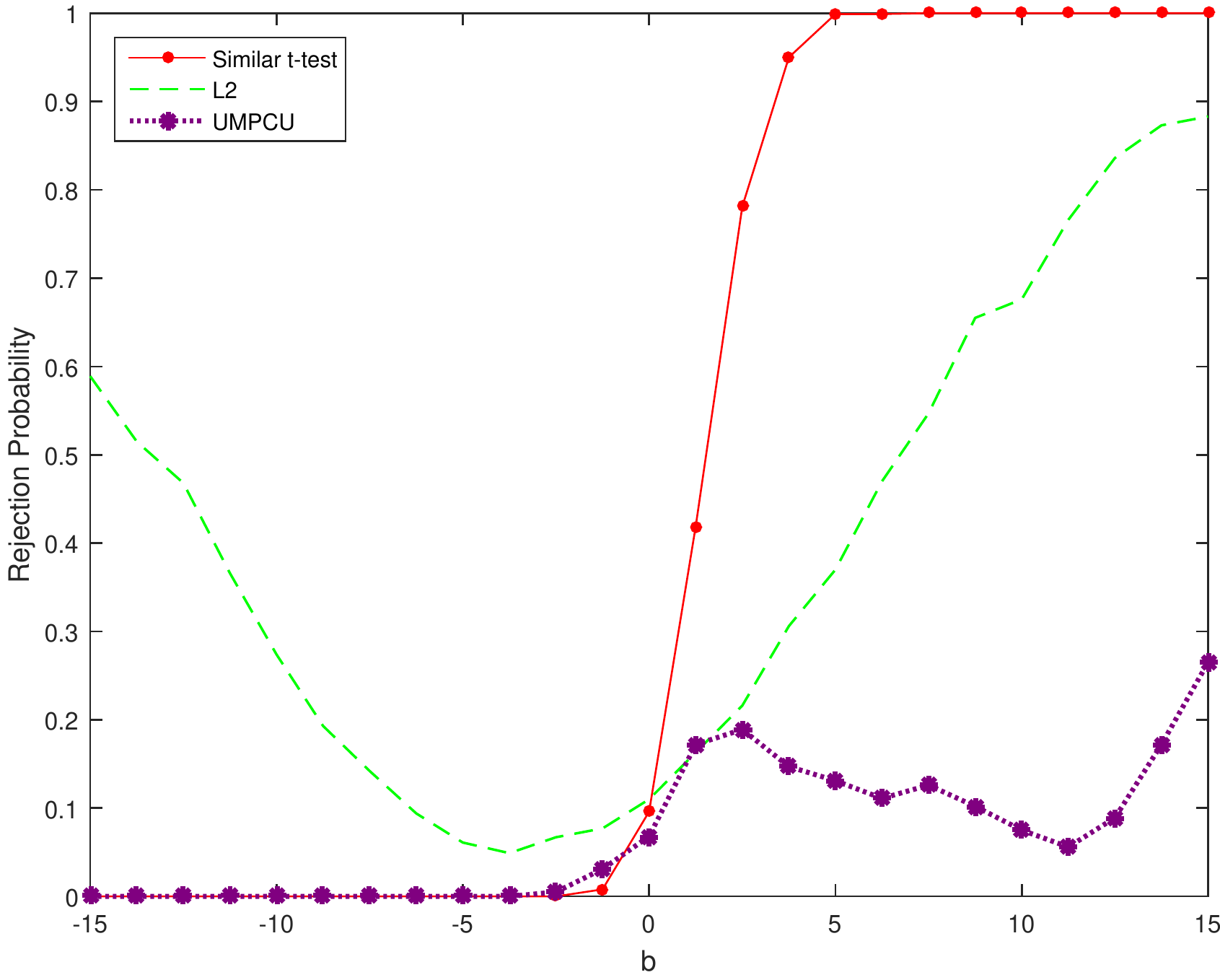} %
\subcaption{$\gamma = 0.5$} \endminipage\hfill \minipage{0.5\textwidth} %
\centering %
\includegraphics[width=6.8cm]{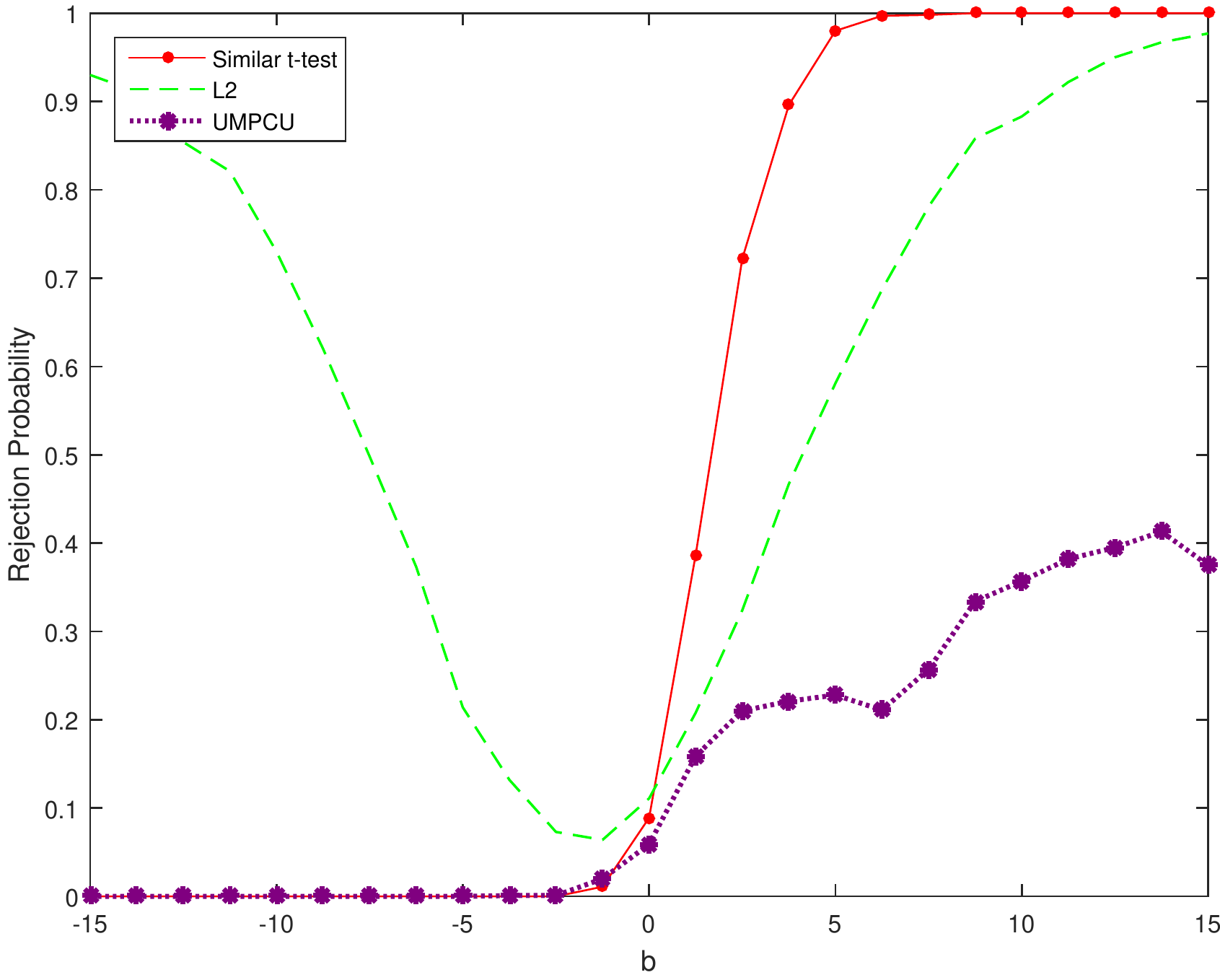} %
\subcaption{$\gamma = 0.85$} \endminipage \hfill \minipage{0.5\textwidth} %
\centering %
\includegraphics[width=6.8cm]{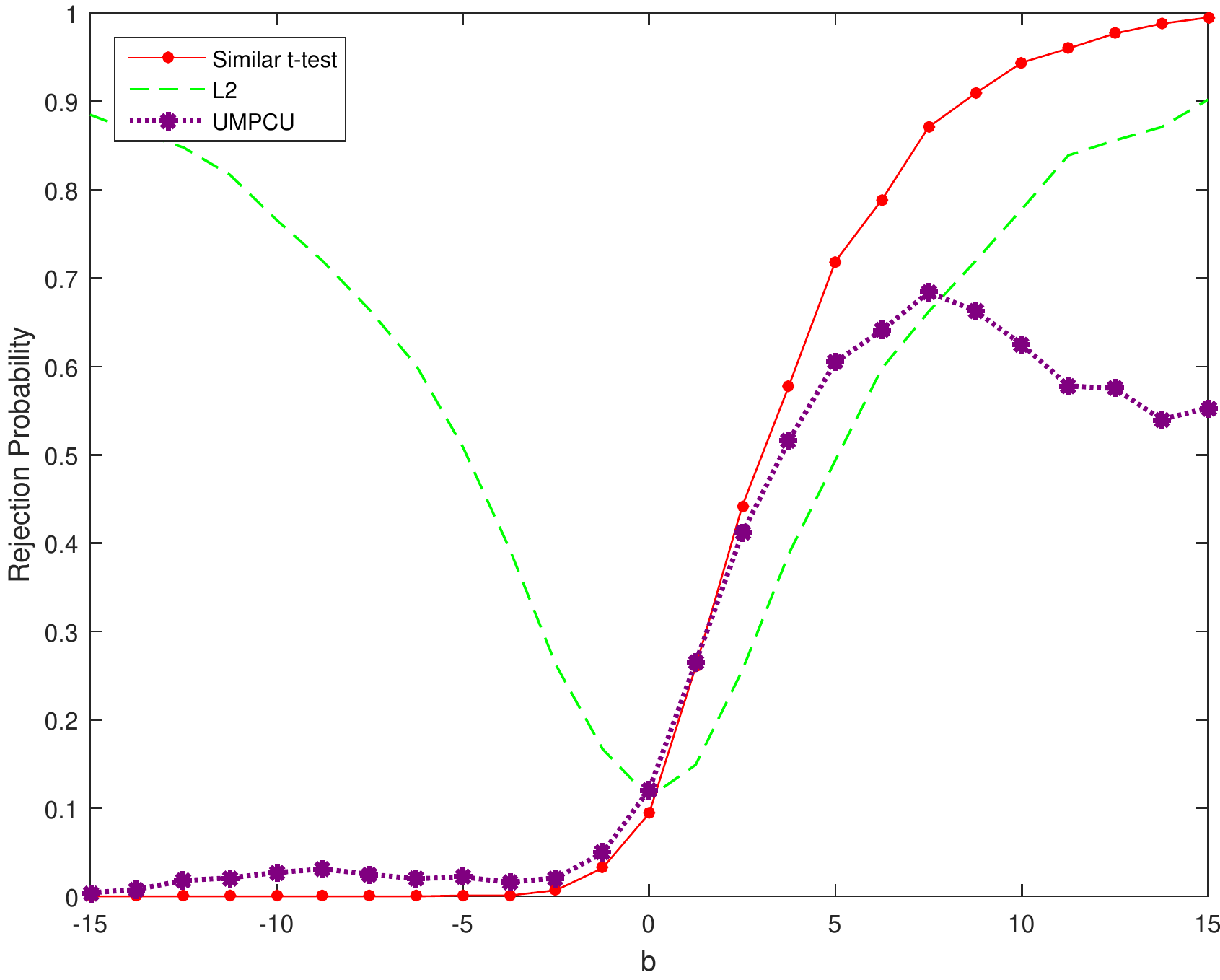} %
\subcaption{$\gamma = 1.0$} \endminipage \hfill \minipage{0.5\textwidth} %
\centering %
\includegraphics[width=6.8cm]{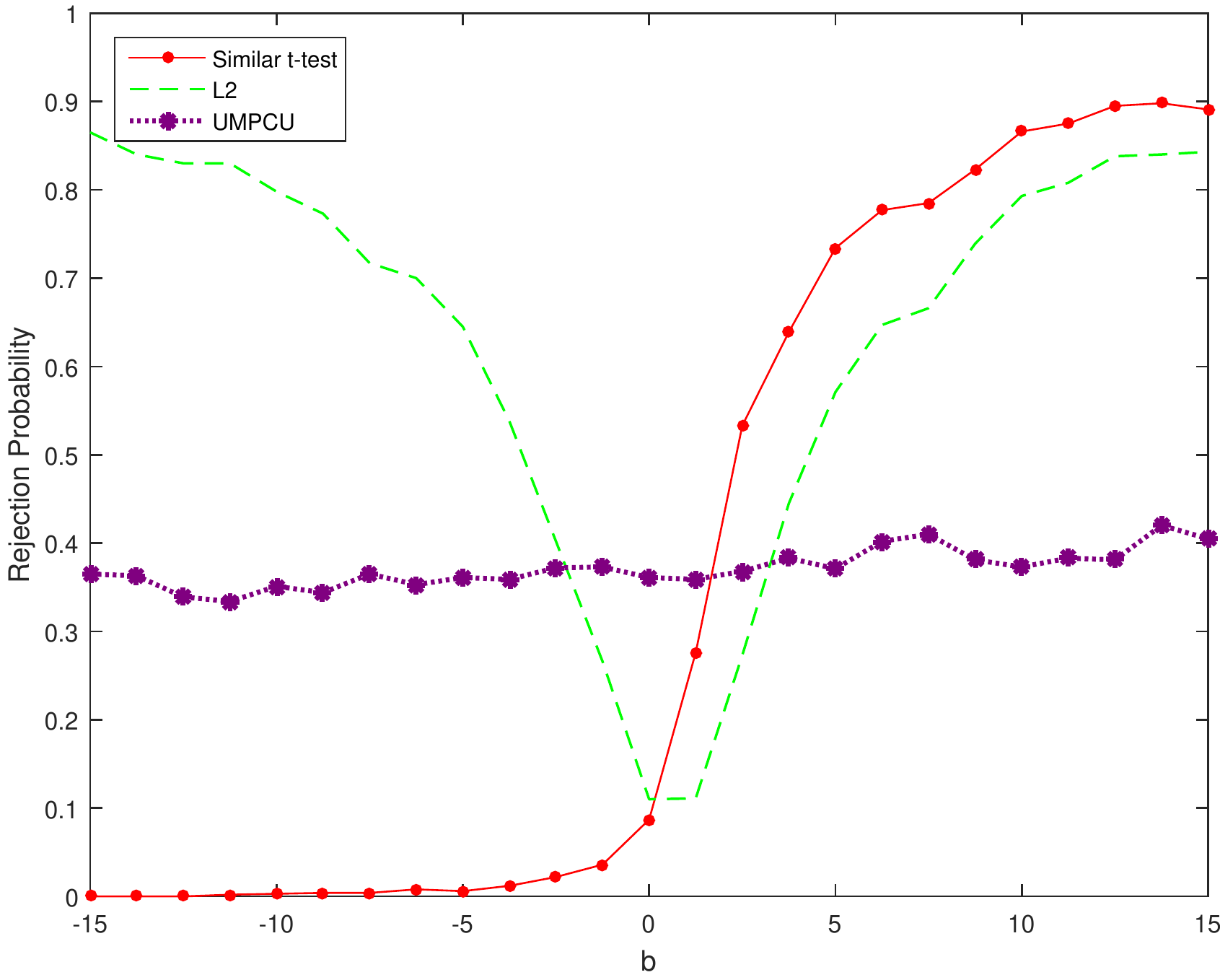} %
\subcaption{$\gamma = 1.05$} \endminipage\hfill \minipage{0.5\textwidth} %
\centering %
\includegraphics[width=6.8cm]{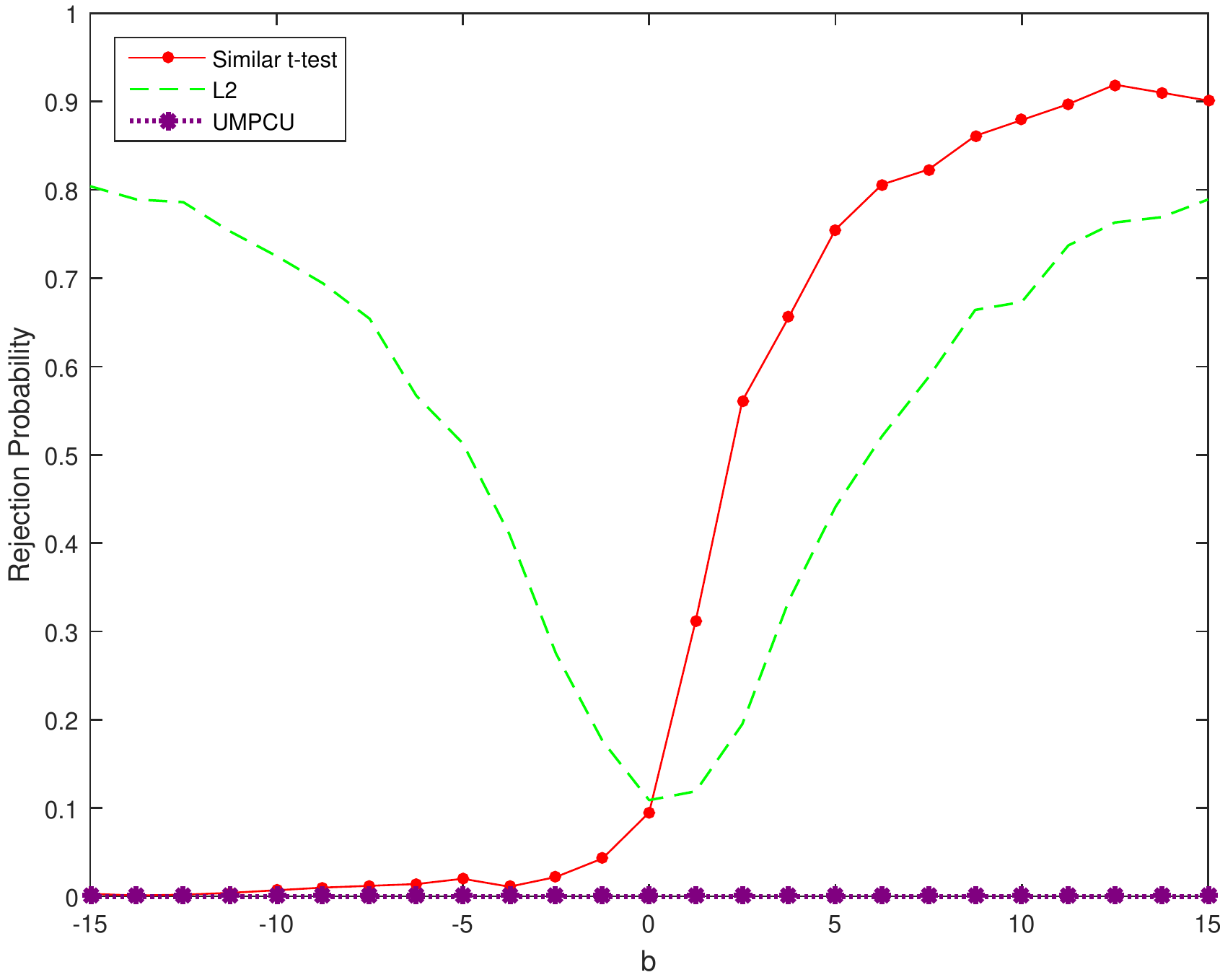} %
\subcaption{$\gamma = 1.3$} \endminipage\hfill
\end{figure}

\begin{figure}[tbh]
\caption{Power ($\protect\rho =-0.5$)}
\label{fig:Power -0.5 suppl}
\bigskip \minipage{0.5\textwidth} \centering %
\includegraphics[width=6.8cm]{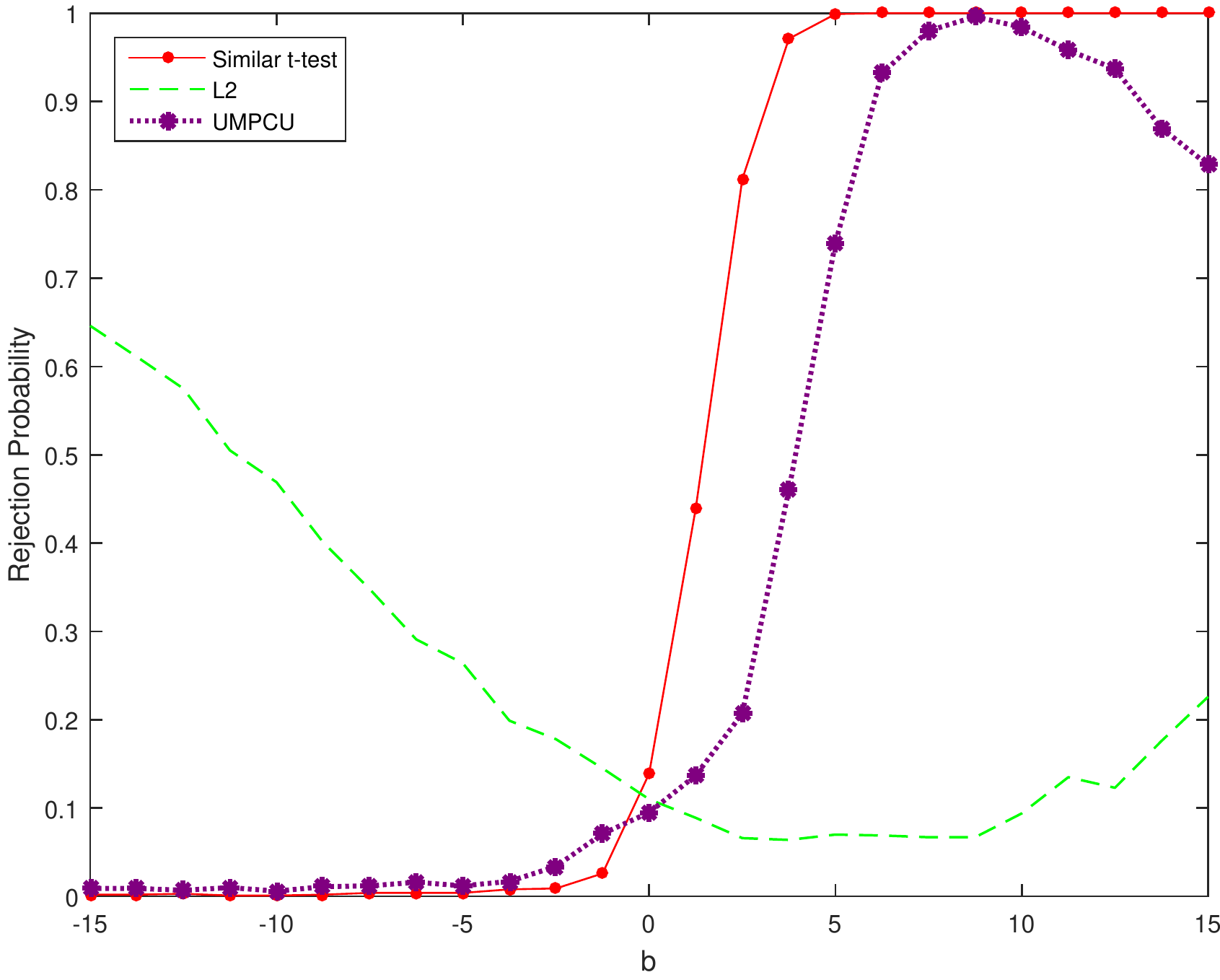} %
\subcaption{$\gamma = 0$} \endminipage \hfill \minipage{0.5\textwidth} %
\centering %
\includegraphics[width=6.8cm]{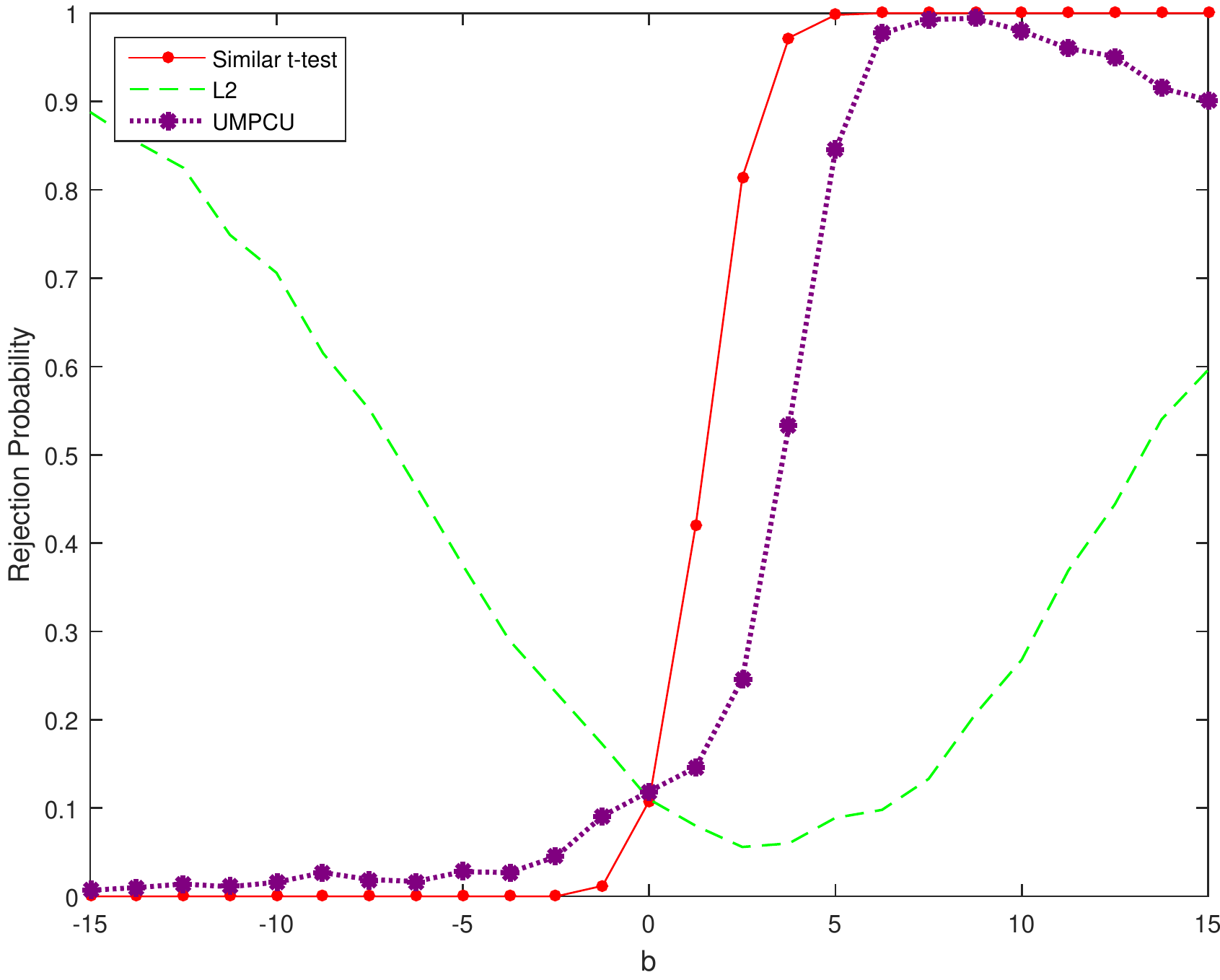} %
\subcaption{$\gamma = 0.5$} \endminipage\hfill \minipage{0.5\textwidth} %
\centering %
\includegraphics[width=6.8cm]{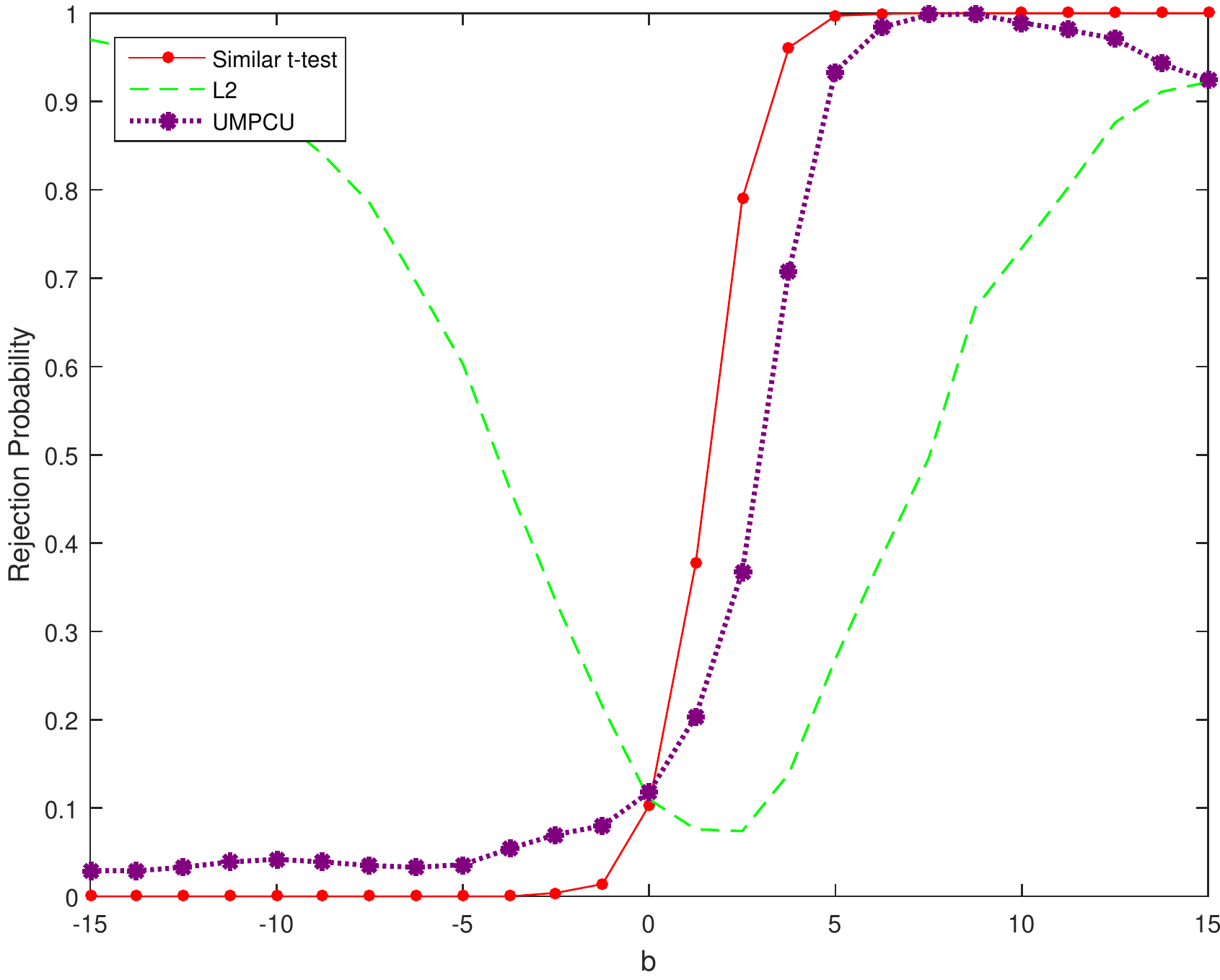} %
\subcaption{$\gamma = 0.85$} \endminipage \hfill \minipage{0.5\textwidth} %
\centering %
\includegraphics[width=6.8cm]{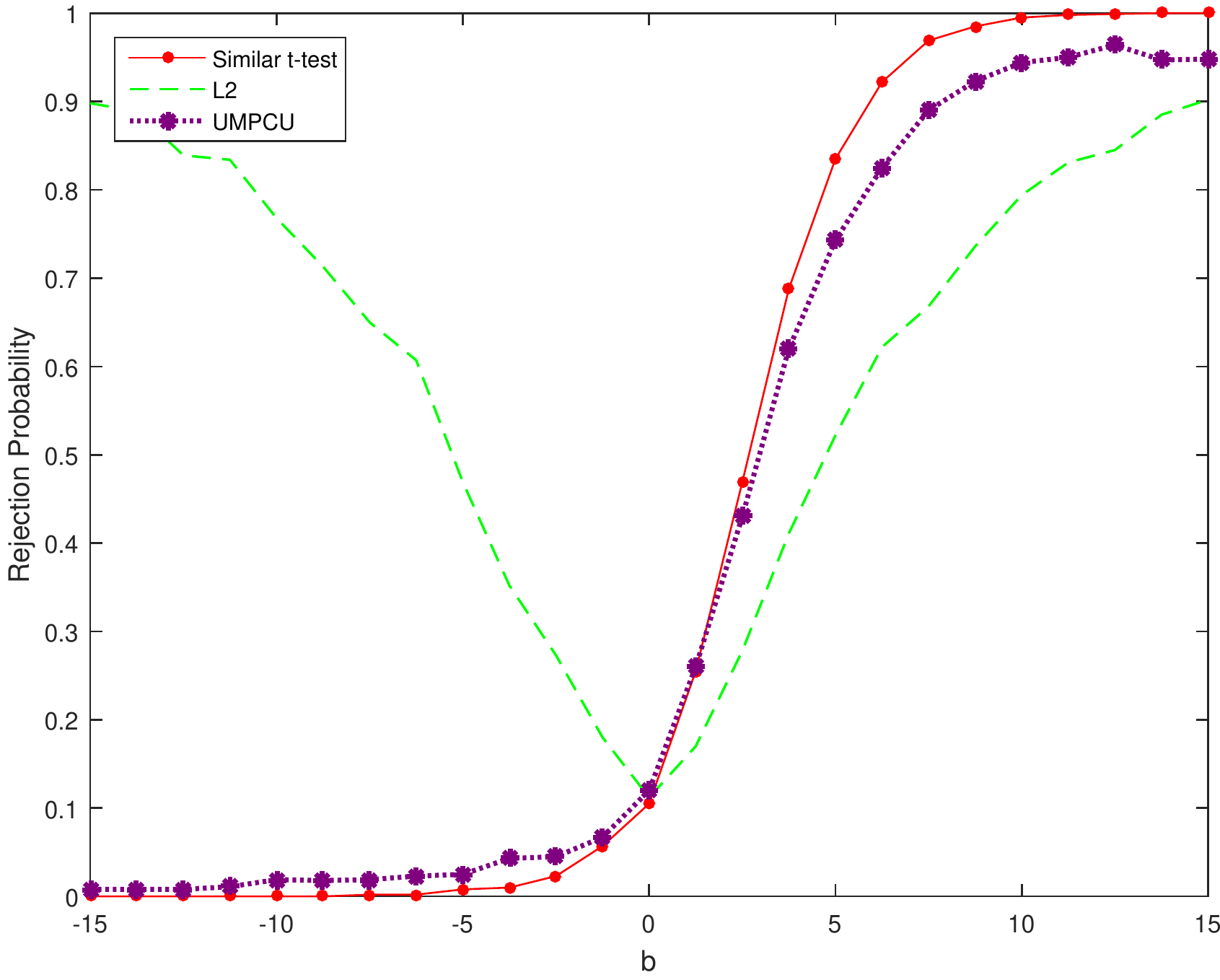} %
\subcaption{$\gamma = 1.0$} \endminipage \hfill \minipage{0.5\textwidth} %
\centering %
\includegraphics[width=6.8cm]{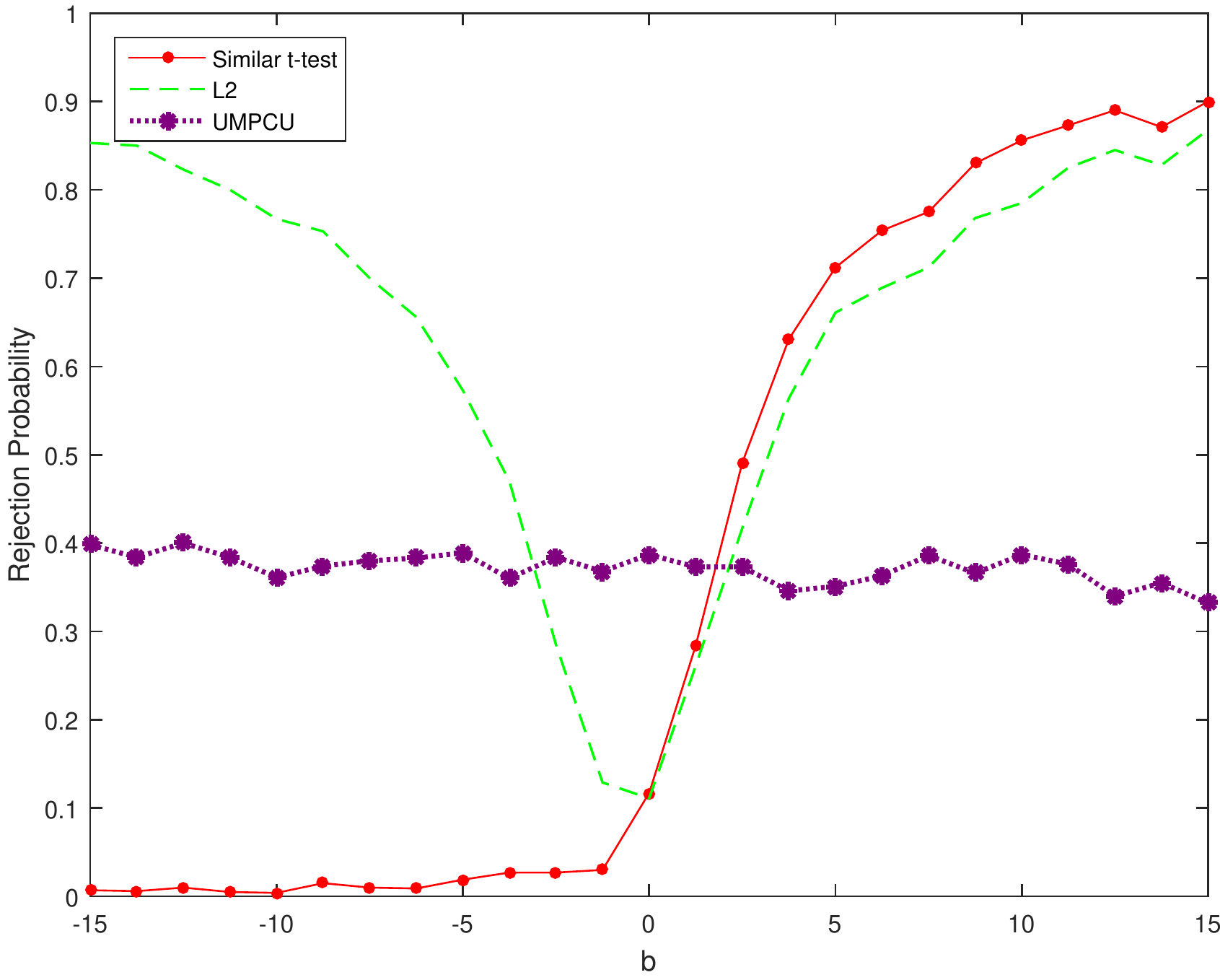} %
\subcaption{$\gamma = 1.05$} \endminipage\hfill \minipage{0.5\textwidth} %
\centering %
\includegraphics[width=6.8cm]{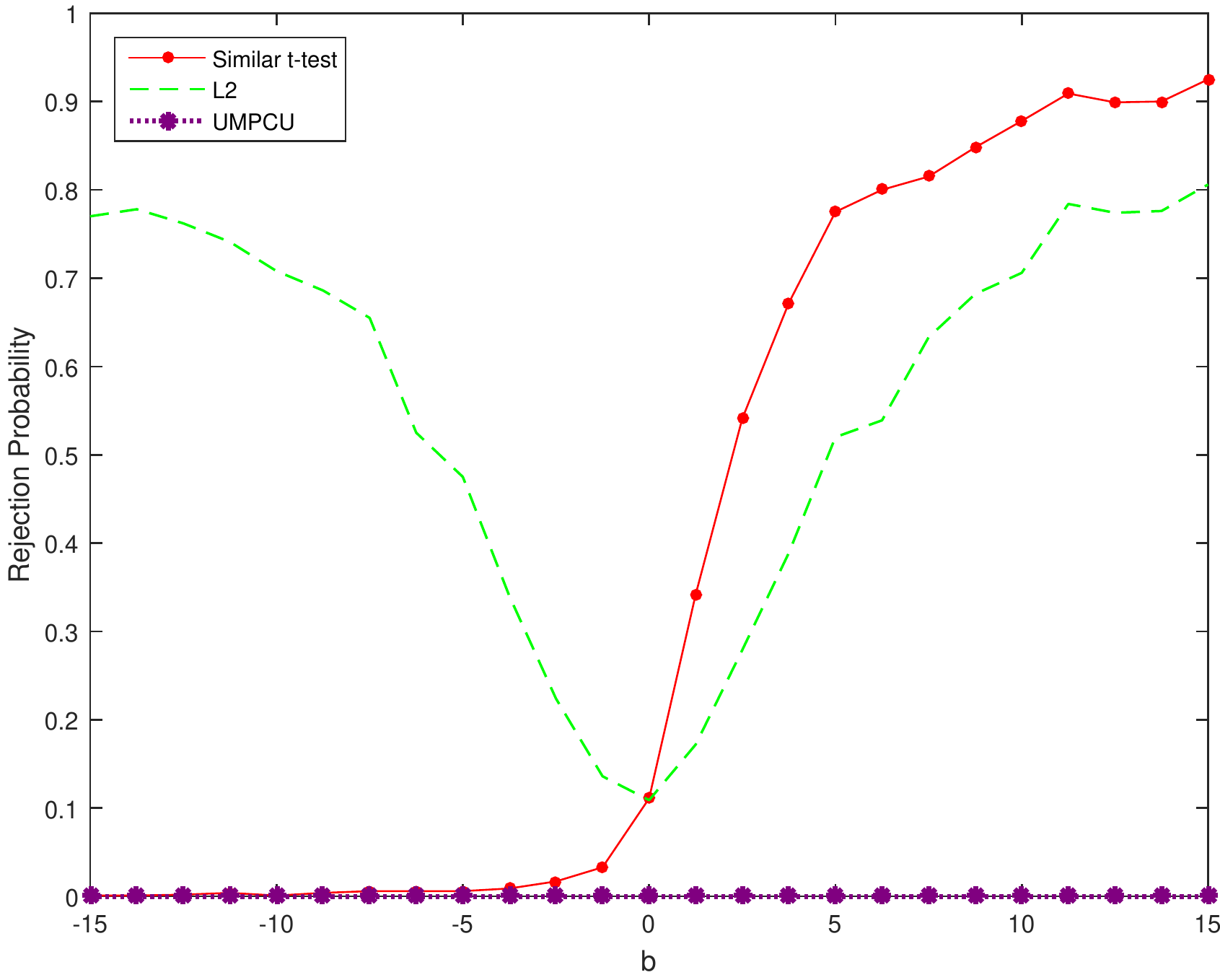} %
\subcaption{$\gamma = 1.3$} \endminipage\hfill
\end{figure}

\begin{figure}[tbh]
\caption{Power ($\protect\rho =0.95$)}
\label{fig:Power 0.95 suppl}
\bigskip \minipage{0.5\textwidth} \centering %
\includegraphics[width=6.8cm]{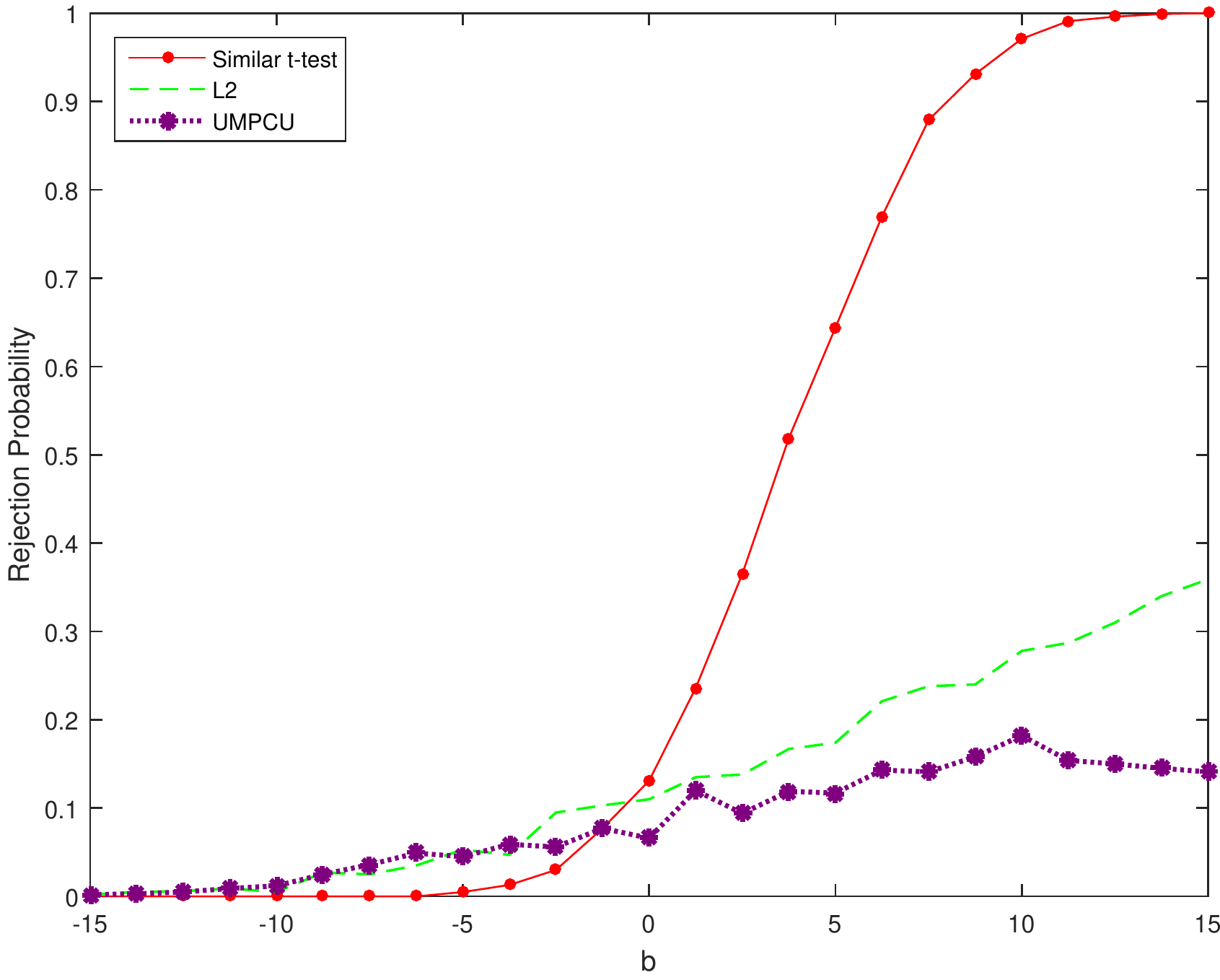} %
\subcaption{$\gamma = 0$} \endminipage \hfill \minipage{0.5\textwidth} %
\centering %
\includegraphics[width=6.8cm]{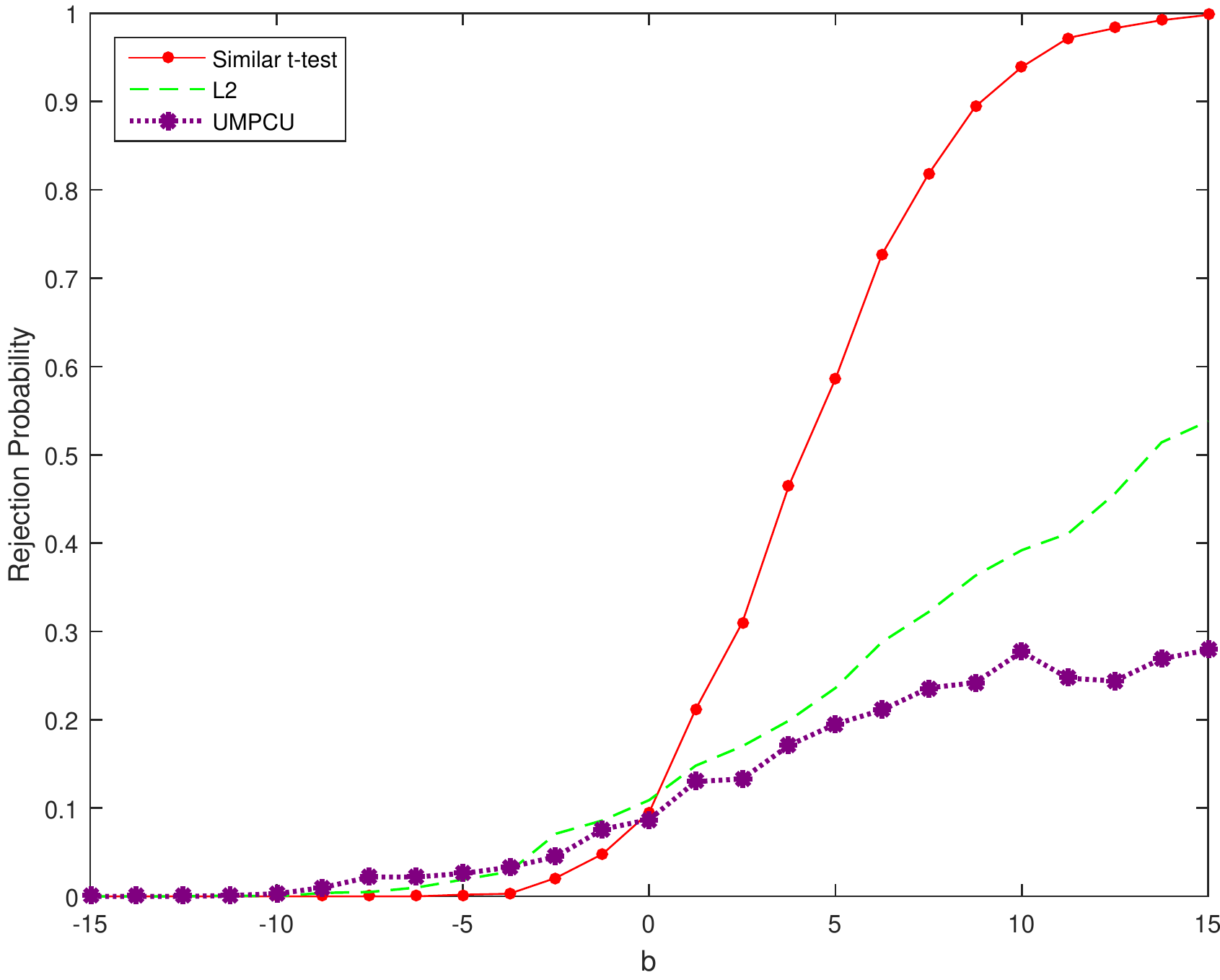} %
\subcaption{$\gamma = 0.5$} \endminipage\hfill \minipage{0.5\textwidth} %
\centering %
\includegraphics[width=6.8cm]{{mm_power_comparison_gamma_0.85_rho_0.95-eps-converted-to}.pdf} %
\subcaption{$\gamma = 0.85$} \endminipage \hfill \minipage{0.5\textwidth} %
\centering %
\includegraphics[width=6.8cm]{{mm_power_comparison_gamma_1_rho_0.95-eps-converted-to}.pdf} %
\subcaption{$\gamma = 1.0$} \endminipage \hfill \minipage{0.5\textwidth} %
\centering %
\includegraphics[width=6.8cm]{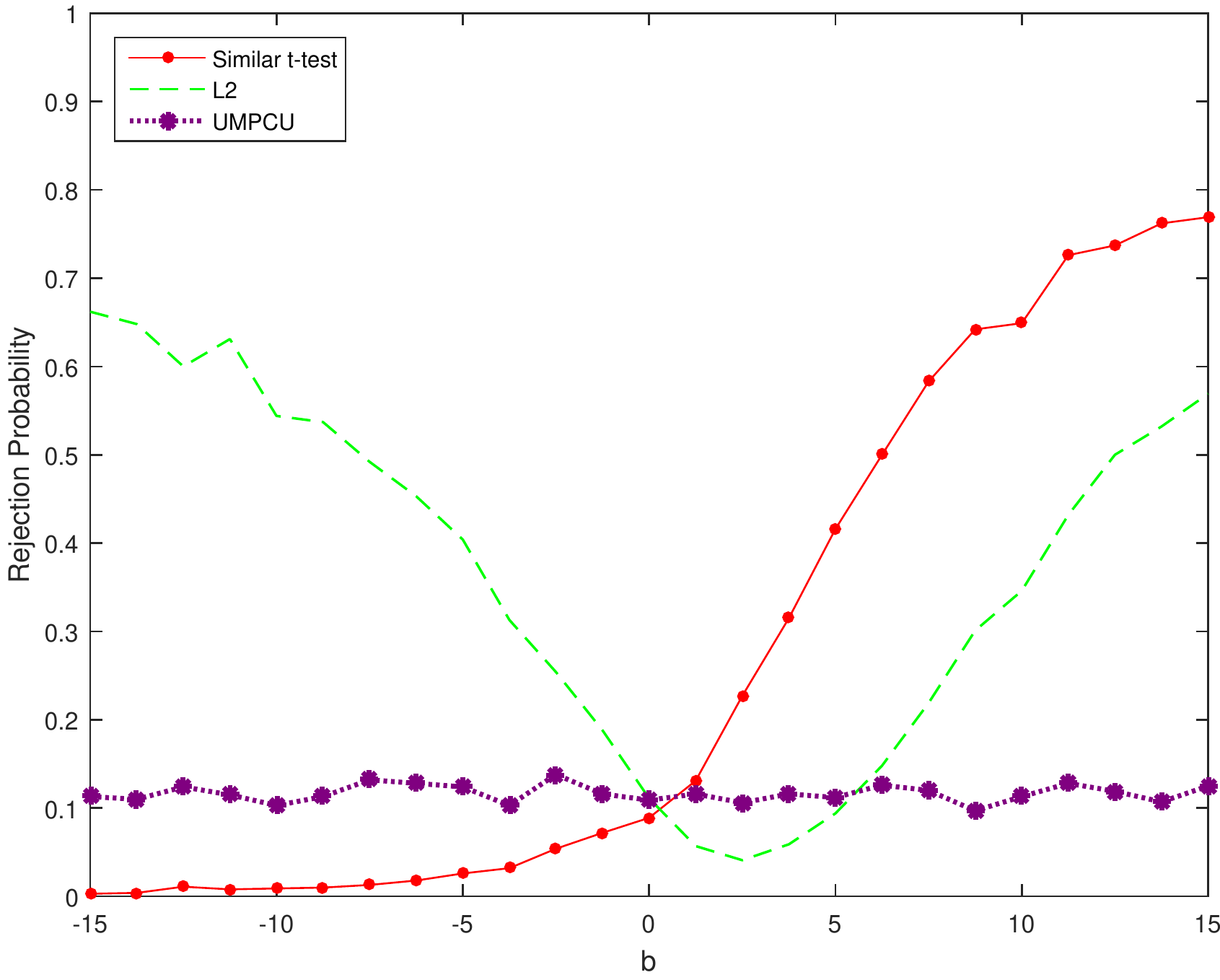} %
\subcaption{$\gamma = 1.05$} \endminipage\hfill \minipage{0.5\textwidth} %
\centering %
\includegraphics[width=6.8cm]{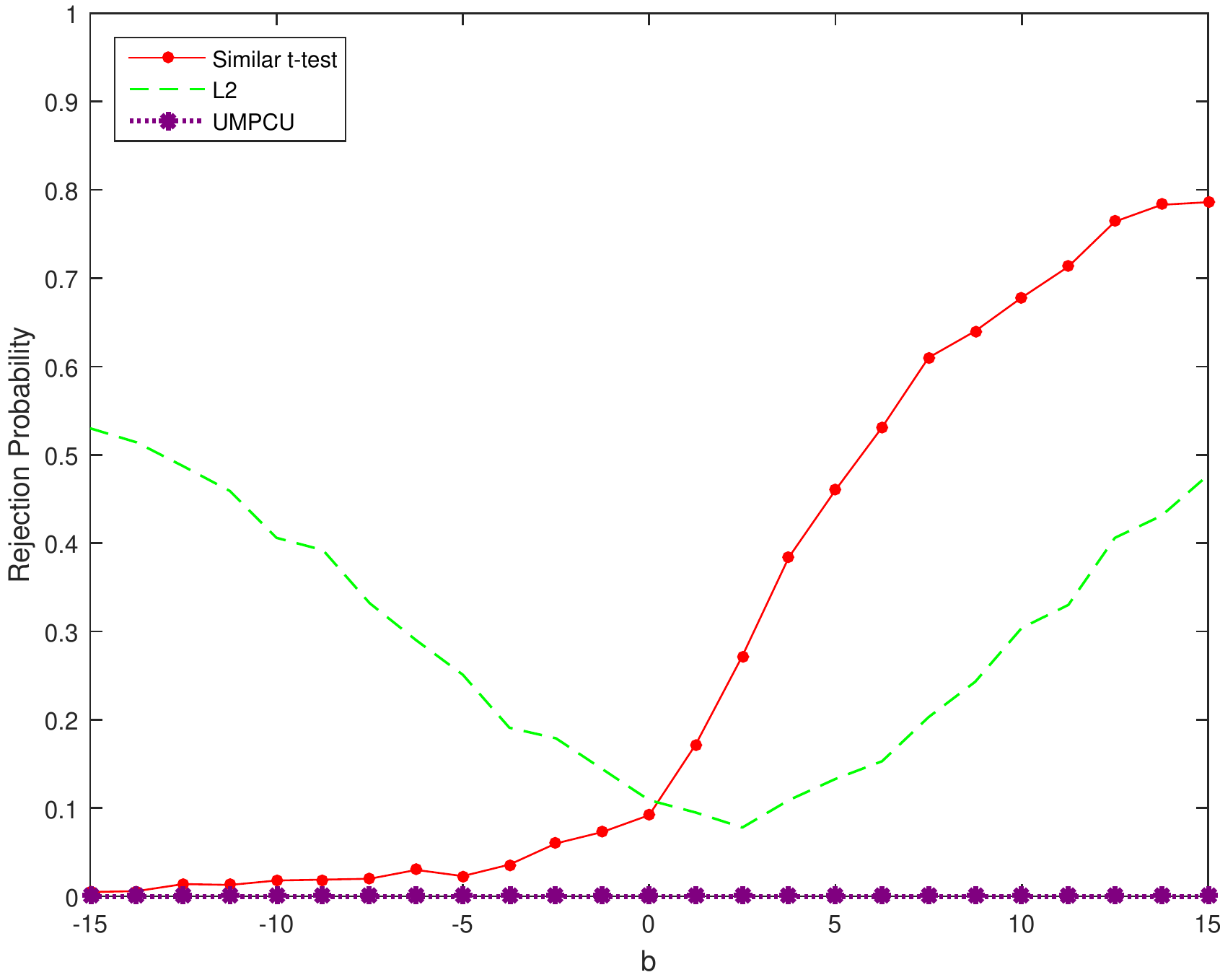} %
\subcaption{$\gamma = 1.3$} \endminipage\hfill
\end{figure}

\begin{figure}[tbh]
\caption{Power ($\protect\rho =-0.95$)}
\label{fig:Power -0.95 suppl}
\bigskip \minipage{0.5\textwidth} \centering %
\includegraphics[width=6.8cm]{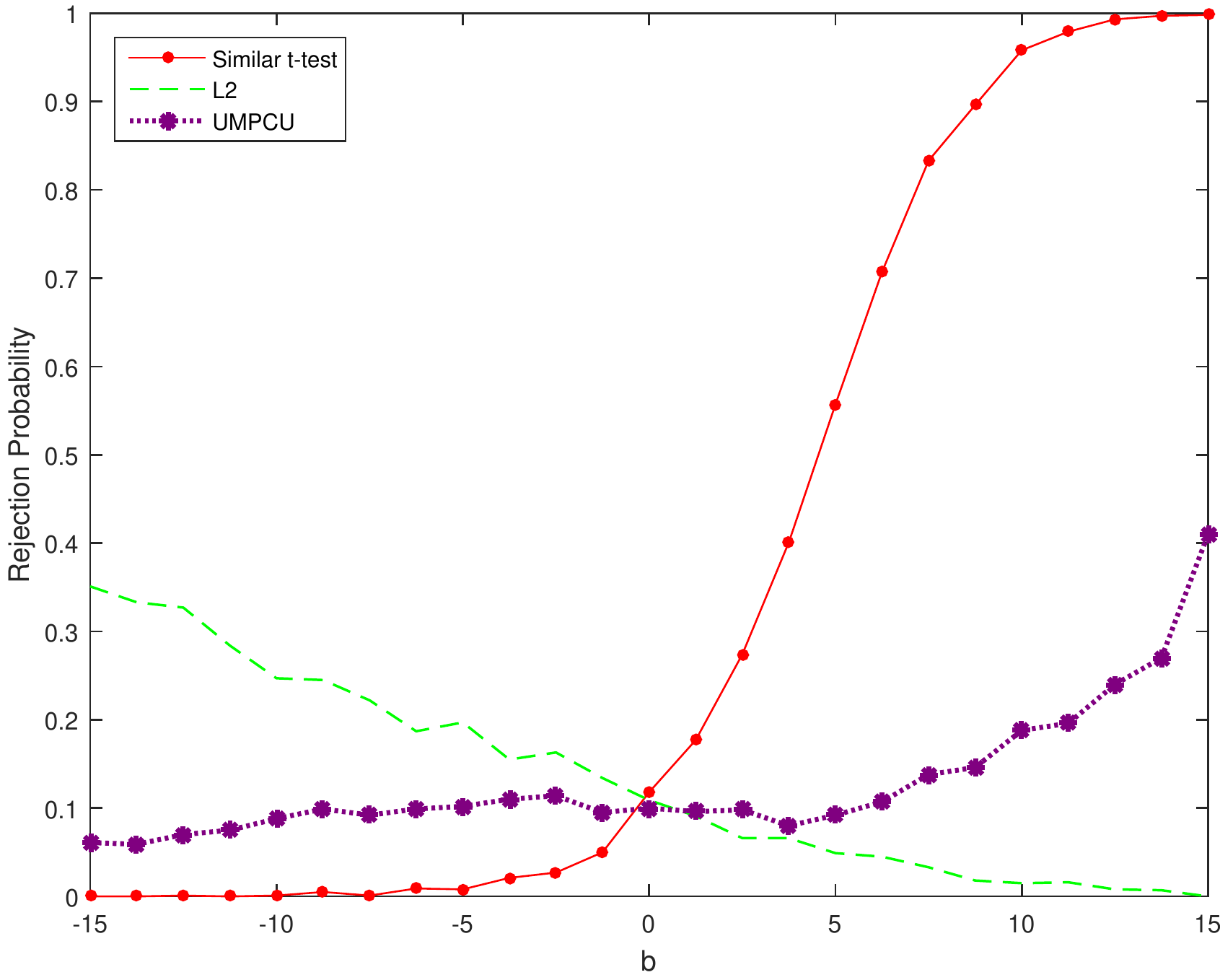} %
\subcaption{$\gamma = 0$} \endminipage \hfill \minipage{0.5\textwidth} %
\centering %
\includegraphics[width=6.8cm]{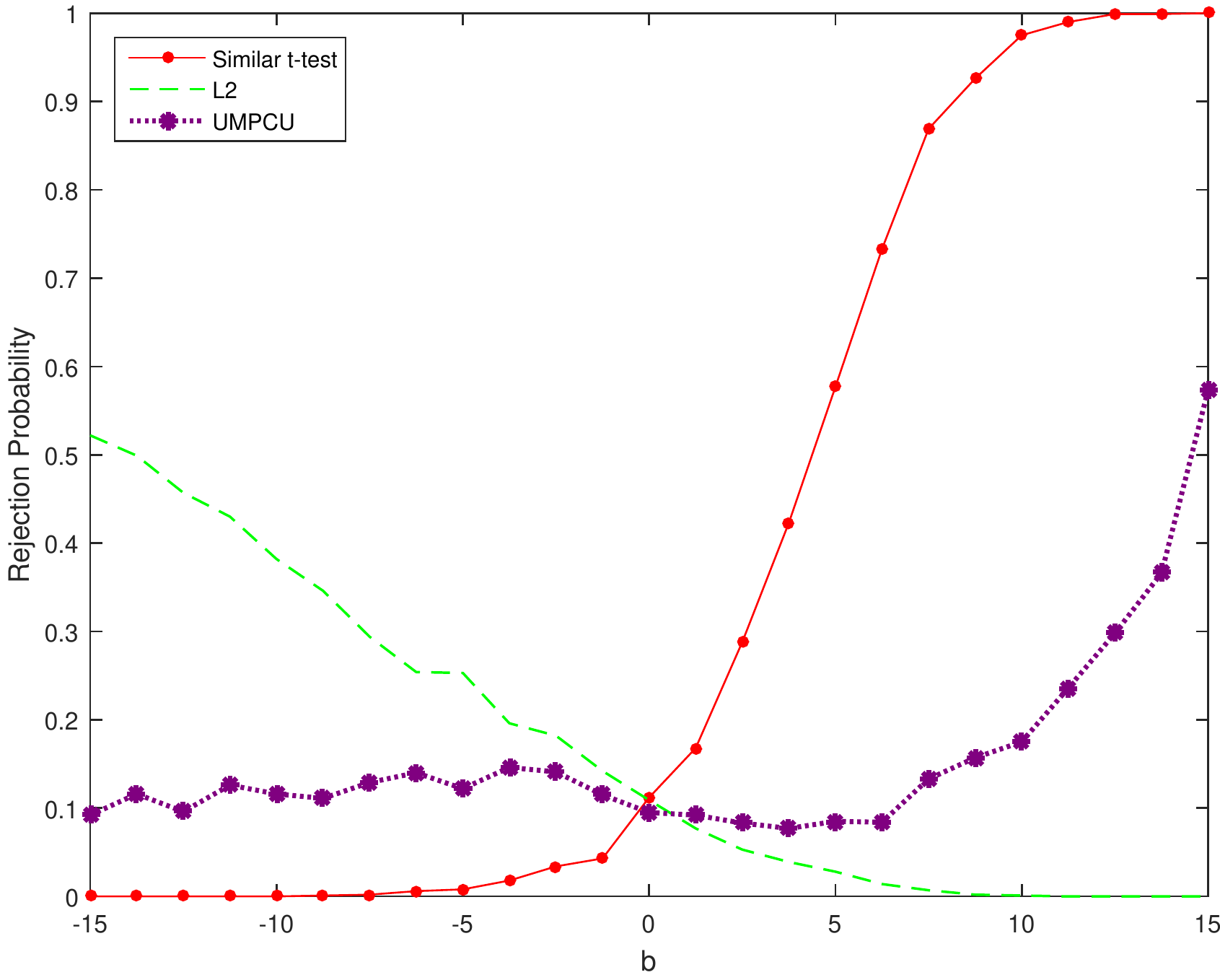} %
\subcaption{$\gamma = 0.5$} \endminipage\hfill \minipage{0.5\textwidth} %
\centering %
\includegraphics[width=6.8cm]{{mm_power_comparison_gamma_0.85_rho_-0.95-eps-converted-to}.pdf}
\subcaption{$\gamma = 0.85$} \endminipage \hfill \minipage{0.5\textwidth} %
\centering %
\includegraphics[width=6.8cm]{{mm_power_comparison_gamma_1_rho_-0.95-eps-converted-to}.pdf} %
\subcaption{$\gamma = 1.0$} \endminipage \hfill \minipage{0.5\textwidth} %
\centering %
\includegraphics[width=6.8cm]{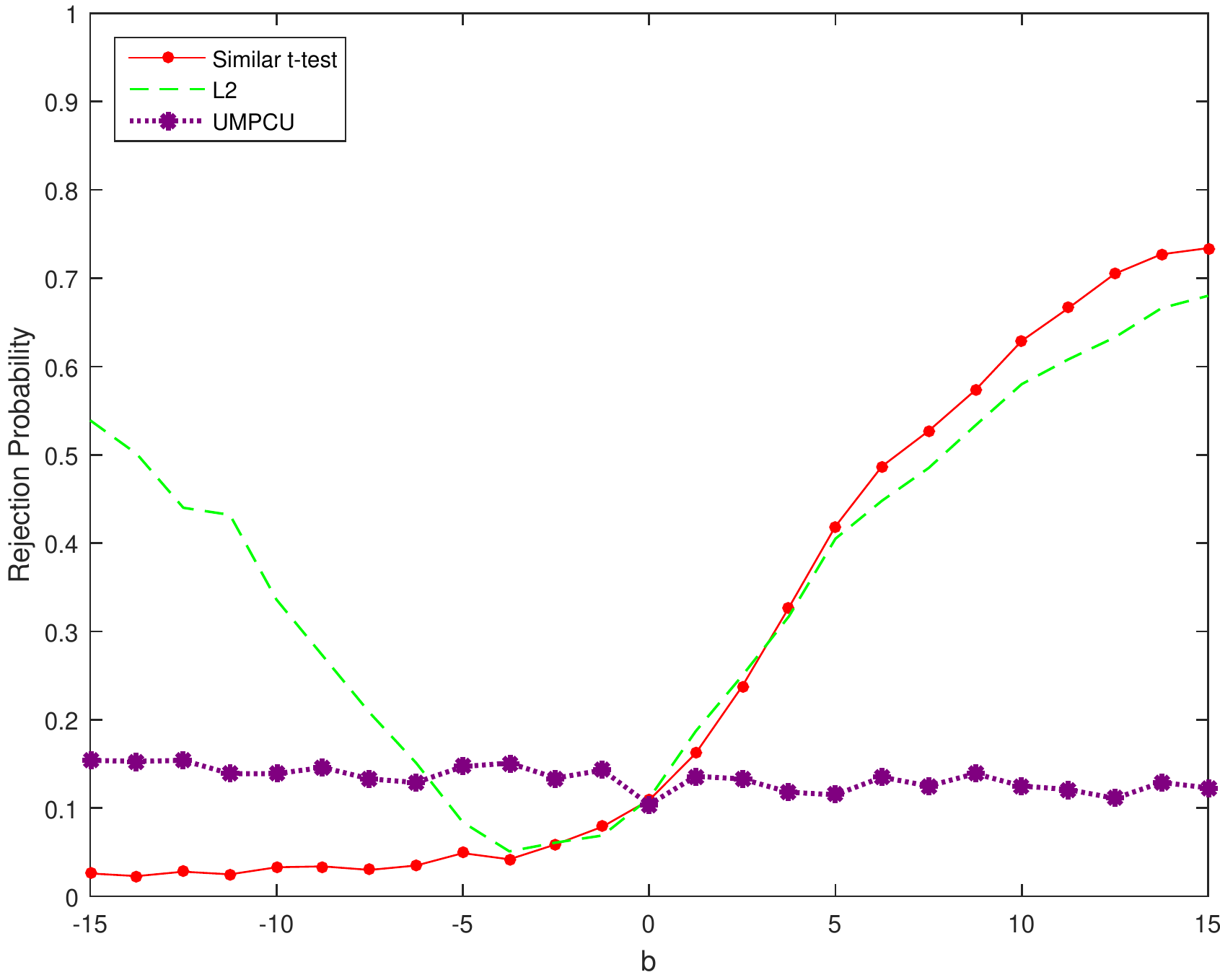}
\subcaption{$\gamma = 1.05$} \endminipage\hfill \minipage{0.5\textwidth} %
\centering %
\includegraphics[width=6.8cm]{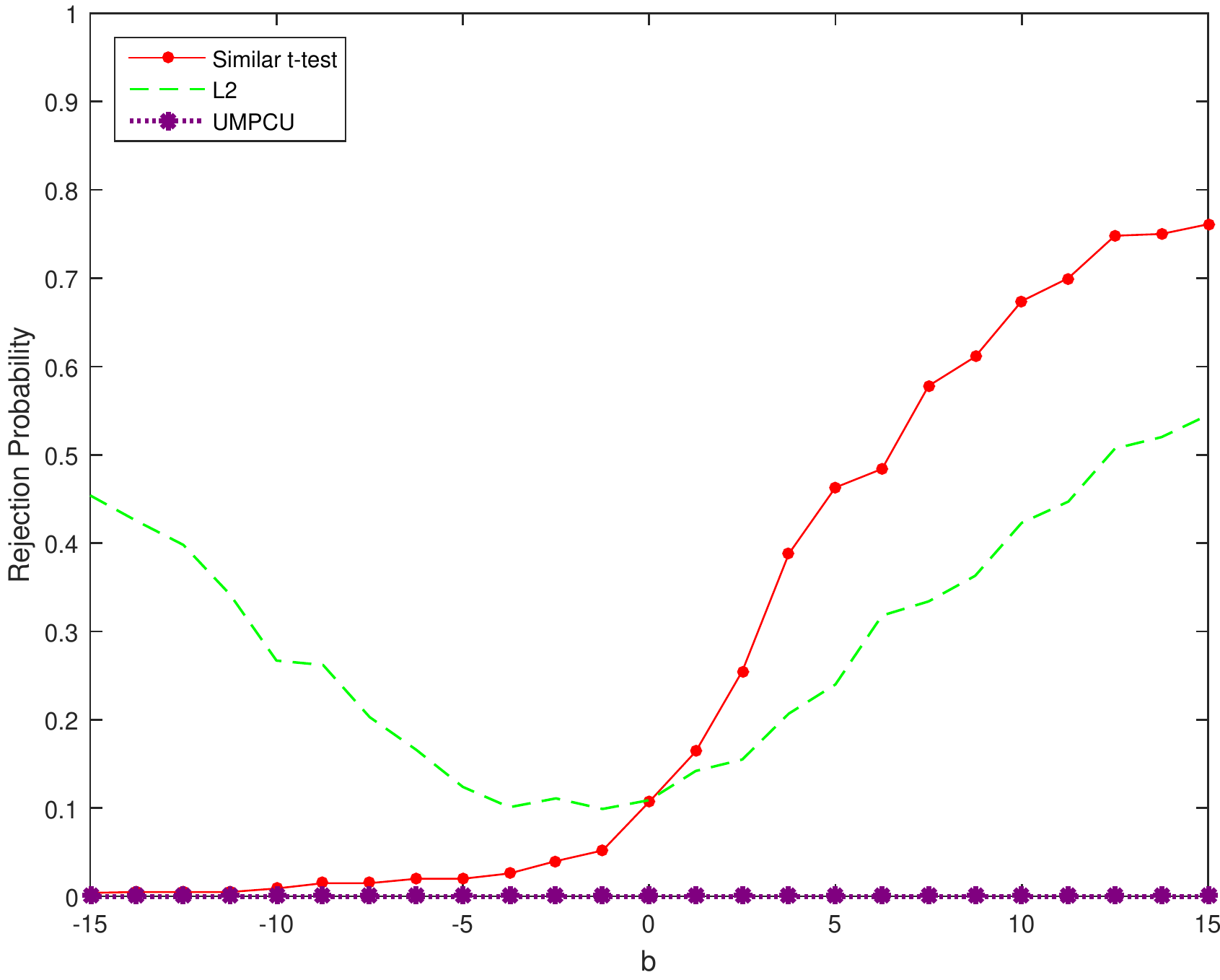} %
\subcaption{$\gamma = 1.3$} \endminipage\hfill
\end{figure}

Figures \ref{fig:Power 0.5 suppl}-\ref{fig:Power -0.95 suppl} present power
for different levels of $\gamma $ and $\rho $. The results further support
the similar t-test in terms of size and power. The behavior of the $L_{2}$
test is sensitive to the values of $c$, sometimes behaving as a one-sided
and other times as a two-sided test. At the value $b=0$, the UMPCU test
presents null rejection probabilities above the nominal level when $c=5$,
and close to zero when $c=30$. This problem is possibly due to the
inaccuracy of the algorithm used by \citet{JanssonMoreira06} for integral
approximations. Although this problem can be conceivably fixed by using
better numerical approximations for $c>0$, the power of the UMPCU test is
considerably below that found by the similar t-test for $c\leq 0$. Overall,
the similar t-test has correct size and better power than the competing
tests.

\bibliographystyle{Chicago}
\bibliography{References}